\documentclass[10pt]{article}

\usepackage{definitions}

\usepackage{etoolbox}
\usepackage[section]{placeins}
\preto\subsection{\FloatBarrier}
\newcommand{\cX}{\mathcal{X}}
\newcommand{\cY}{\mathcal{Y}}
\newcommand{\cZ}{\mathcal{Z}}
\newcommand{\cE}{\mathcal{E}}
\newcommand{\cR}{\mathcal{R}}
\newcommand{\cB}{\mathcal{B}}
\newcommand{\cQ}{\mathcal{Q}}
\newcommand{\cS}{\mathcal{S}}
\newcommand{\cF}{\mathcal{F}}
\newcommand{\Ep}{E_{+}}
\newcommand{\Em}{E_{-}}
\newcommand{\tY}{\tilde{Y}}
\newcommand{\tZ}{\tilde{Z}}
\newcommand{\gp}{\mathrm{gap}}
\newcommand{\diag}{\operatorname{diag}}
\DeclareMathOperator{\dist}{dist}
\DeclareMathOperator{\diam}{diam}

\newcommand{\cnorm}[1]{\norm{#1}_{\Lambda}}

\pgfplotsset{compat=1.18}

\title{Reformulation Techniques for PDLP with Applications \\ to QUBO and Joint Replenishment}
\author{Jacob Dentes \thanks{MIT, Sloan School of Management (dentes@mit.edu).} \and Haihao Lu\thanks{MIT, Sloan School of Management (haihao@mit.edu).}}
\date{}

\begin{document}
\maketitle

\begin{abstract}
Practitioners of linear programming widely believe that reformulating a linear program (LP) offers little benefit to interior-point methods: their iteration counts are famously insensitive to the formulation, so the payoff of reformulation beyond automatic presolve and customized linear solve is modest.
This paper presents evidence that the same is \emph{not} true of first-order solvers built on the primal--dual hybrid gradient method (PDHG), such as PDLP.
Unlike interior-point methods, PDHG-type methods have per-iteration costs governed by the matrix representation of the formulation and iteration counts governed by its geometry, and both change under equivalent reformulations.
We describe a general reformulation recipe: constraints drawn from a catalogue of \emph{projection-friendly} polyhedra can be absorbed into the projection step of PDHG.
We analyze what absorption does to the four constants that govern restarted PDHG: the coupling norm, the termination certificate, the distance scale, and the sharpness modulus.
We perform the reformulation on two case studies.
On the LP relaxation of quadratic unconstrained binary optimization (QUBO), absorption makes any step size below $1$ admissible on every instance, and yields a fully explicit linear rate.
On the LP relaxation of the joint replenishment problem (JRP), whose optimum certifies a $2$-approximation for an NP-hard inventory problem, absorption improves the coupling norm and eliminates every matrix-vector product.
Both reformulations run on GPUs with custom kernels in place of sparse matrix-vector products. We benchmark both against the GPU LP solver cuPDLPx applied to the original formulations. By solve time, the reformulated solvers are $6.4$--$9.6\times$ faster on QUBO and $6.1$--$28.9\times$ faster on the JRP. They take $1.7$--$4.3\times$ and $13.6$--$52.3\times$ fewer iterations respectively.
\end{abstract}

\section{Introduction}\label{sec:intro}

Linear programming is a foundational tool of operations research, and for decades the simplex and interior-point methods have dominated the space of LP solvers \cite{bixby_solving_2002}.
Both rely on matrix factorizations whose memory traffic and fill-in scale poorly on modern massively parallel hardware.
First-order methods, which trade factorizations for matrix-vector products, have recently emerged as a practical alternative for large-scale LP \cite{lu_first-order_2024,nesterov_lectures_2018}.
PDLP, a development of PDHG \cite{chambolle_first-order_2011}, demonstrated that a carefully engineered first-order method is competitive on broad LP benchmarks \cite{applegate_pdlp_2026}.
Its computational core is matrix-vector multiplication, so PDHG suits GPUs.
cuPDLP.jl \cite{lu_cupdlpjl_2024} removed the remaining sequential components of PDLP and showed that GPU implementations can achieve numerical performance on par with mature CPU solvers, while cuPDLP-C \cite{lu_cupdlp-c_2024} provided a mature C implementation of this approach. D-PDLP \cite{li_dpdlp_2026} scaled PDLP to multi-GPU systems, while HPR-LP \cite{chen_hpr-lp_2026} pursued an alternative first-order approach based on the Halpern Peaceman--Rachford method.
PDHG-type LP solvers now ship in production codes, including Gurobi's PDHG \cite{gurobi_optimization_llc_gurobi_2026}, NVIDIA's cuOpt \cite{nvidia_corporation_nvidia_2026}, COPT's PDLP \cite{ge_cardinal_2026}, HiGHS's bundled CuPDLP-C \cite{lu_cupdlp-c_2024,highs_development_team_highs_2024}, FICO Xpress's PDHG \cite{fair_isaac_corporation_fico_2024}, and Knitro's PDLP \cite{artelys_algorithms_2026}.

This paper presents a lever that the LP community has largely written off: \emph{reformulation}.
Folklore among practitioners holds that reformulating an LP offers limited benefit to interior-point methods.
Their iteration counts are remarkably stable, nearly independent of instance size or data.
Formulation affects performance mainly through the sparsity of the factorized systems.
Customized linear solve and automatic presolve already capture most of that benefit \cite{gondzio_interior_2012,wright_primal-dual_1997,bixby_solving_2002}.
Under this folklore there is little reason for a modeler to hand-craft an equivalent formulation.

PDHG-type methods are different in two ways.
First, their per-iteration cost \emph{is} the cost of multiplying by the constraint matrix, so a reformulation that shrinks or restructures that matrix changes the price of every iteration directly.
Second, their convergence speed is governed by the \emph{geometry} of the formulation.
Restarted PDHG converges linearly at a rate set by sharpness constants of the
primal--dual formulation \cite{applegate_faster_2023,lu_restarted_2024}.
Those constants are closely tied to Hoffman constants and level-set geometry
\cite{xiong_computational_2026,xiong_role_2024}, and equivalent reformulations change
them.
Reformulation therefore has two levers on a PDHG-type solver where it has almost none on an interior-point method.

We evaluate reformulation by absorbing constraints into projected sets.
PDLP projects onto boxes; PDHG needs only closed convex feasible sets with computable projections.
So any subset of the constraints whose induced feasible set has an efficient Euclidean projection can be \emph{absorbed}, moved out of the constraint matrix and into the projection step, instead of dualized.
The solver then runs on a different saddle-point problem for the same LP, with a smaller coupling matrix, fewer dual variables, and different geometry.

Our contributions are summarized as follows.

\begin{itemize}
\item \textbf{A systematic framework for absorption.}
We formalize the absorption procedure and identify classes of polyhedra that admit efficient projections (Section~\ref{sec:reform}).
We also develop a construction for problems involving auxiliary variables that select maxima over groups (Section~\ref{sec:maxselection}), which serves as the basis for our applications in Section~\ref{sec:applications}.

\item \textbf{Algorithmic improvements and convergence guarantees.}
We analyze how absorption affects four key quantities: the coupling norm, the termination certificate, the distance scale, and the sharpness modulus, as well as its potential to reduce per-iteration cost (Section~\ref{sec:why}).
We establish that sharpness holds automatically for polyhedral absorbed formulations and show how a lattice structure in the vertex values of both value polytopes yields an explicit sharpness bound (Corollary~\ref{cor:lattice}), enabling quantitative convergence guarantees.

\item \textbf{Two case studies.}
We apply the framework to the QUBO relaxation (Section~\ref{sec:qubo}) and the joint replenishment LP relaxation (Section~\ref{sec:jrp}), and implement each reformulation as a GPU solver: PRAQ for QUBO and PRIA for joint replenishment.
For QUBO, half-integrality on both sides yields explicit constants and a linear convergence rate (Theorem~\ref{thm:qubo-rate}).
For joint replenishment we establish a sublinear rate whose only data dependence is the instance size (Theorem~\ref{thm:jrp-rate}).
We further evaluate both reformulations against cuPDLPx and Gurobi on the original formulations.
\end{itemize}

\paragraph{Code.}
PRAQ, PRIA, the baselines they are measured against and the benchmark harness are
at \url{https://github.com/MIT-Lu-Lab/reformulation}.

\paragraph{Notation.}
For a positive integer $n$, $[n] := \{1,\dots,n\}$.
For a closed convex set $S$, $\proj_S$ is the projection onto $S$.
For a symmetric positive semidefinite $M$ we write $\norm{v}_M := \sqrt{v^\top M v}$, and
$\diam_M(S)$ and $\dist_M(\cdot,S)$ for the diameter of $S$ and the distance to it in
$\norm{\cdot}_M$.
$\vc(X)$ is the vectorization of a tensor $X$ ordered lexicographically, and $\langle\cdot,\cdot\rangle$ the corresponding inner product.
$\mathbf 1$ is the all-ones vector and $\Delta^r := \{v \ge 0 : \mathbf 1^\top v = r\}$ the simplex of total $r > 0$.

\section{Reformulation for PDHG-type methods}\label{sec:reform}

\subsection{PDHG for linear programming}\label{sec:pdhg-template}
PDHG solves convex--concave saddle-point problems
\begin{equation}\label{eq:saddle}
    \min_{x \in \cX}\ \max_{y \in \cY}\ L(x,y) := \langle c, x\rangle + \langle y,\ q - Ax\rangle,
\end{equation}
with $\cX \subseteq \RR^n$ and $\cY \subseteq \RR^m$ closed convex, by alternating proximal gradient steps \cite{chambolle_first-order_2011} under positive diagonal scalings $\Lambda_{\cX},\Lambda_{\cY}$:
\begin{equation}\label{eq:pdhg}
    x^{k+1} := \proj_{\cX}\!\pran{x^k - \tau\Lambda_{\cX}^{-1}\pran{c - A^\top y^k}},
    \qquad
    y^{k+1} := \proj_{\cY}\!\pran{y^k + \sigma\Lambda_{\cY}^{-1}\pran{q - A(2x^{k+1} - x^k)}},
\end{equation}
which converges whenever $\tau\sigma\cnorm{A}^2 < 1$ for the \emph{effective coupling norm}
\begin{equation}\label{eq:cnorm}
    \cnorm{A} \;:=\; \norm{\Lambda_{\cY}^{-1/2}A\Lambda_{\cX}^{-1/2}}_2 .
\end{equation}
Both projections are taken in $\norm{\cdot}_{\Lambda_{\cX}}$ and $\norm{\cdot}_{\Lambda_{\cY}}$, which makes one application of \eqref{eq:pdhg} a proximal-point step \cite{he_convergence_2012}.
$\Lambda_{\cX}$ rescales the primal step, which moves $\cnorm{A}$ and the admissible step size. At $\Lambda_{\cX} = \Lambda_{\cY} = I$ this is the unpreconditioned method and $\cnorm{A} = \norm{A}_2$.
We write $\tau = \eta/\omega$ and $\sigma = \eta\omega$, splitting the steps into a \emph{step size} $\eta$ and a \emph{primal weight} $\omega$, and denote one application of \eqref{eq:pdhg} by $T_\omega$.

A general LP fits \eqref{eq:saddle} by taking $\cX$ to be the variable bounds (a box), $\cY$ a product of half-lines dual to the rows, and $A$ the full constraint matrix.
This is the \emph{standard form}, the formulation solved by PDLP and its successors \cite{applegate_pdlp_2026,lu_cupdlpx_2025}, and each of its iterations costs two sparse matrix-vector products plus two trivial projections.

\subsection{Absorbing constraints into the projection}\label{sec:absorb}

PDLP projects onto boxes, but nothing in \eqref{eq:pdhg} requires that.
Beyond the closed convexity of $\cX$ and $\cY$, PDHG needs only that $\proj_{\cX}$ and $\proj_{\cY}$ be computable.

Consider an LP with the constraint system
\begin{equation}\label{eq:lp-def}
    \min c^\top x \qquad \mathrm{s.t.}\ A x \geq q, \ x \in \text{box}
\end{equation}
A row of \eqref{eq:lp-def} can be \emph{dualized}, in which case PDHG carries a multiplier for it, or \emph{absorbed}, in which case $\proj_{\cX}$ maintains it and no multiplier exists.

\begin{mydef}[Absorption split]\label{def:split}
An \emph{absorption split} of \eqref{eq:lp-def} is a partition of its rows into an absorbed set $\cE$ and a dualized set $\cR$ such that
    \[
\cX := \{x \in \mathrm{box} : A_{\cE}x \ge q_{\cE}\}
    \]
has a computable projection.
\end{mydef}

Under the split, \eqref{eq:lp-def} reads
\begin{equation}\label{eq:lp-split}
    \min c^\top x \qquad \mathrm{s.t.}\  A_{\cE}x \ge q_{\cE}, \ A_{\cR}x \ge q_{\cR}, \ x \in \text{box},
\end{equation}
where $\proj_{\cX}$ maintains the first group of rows and the second keeps its multipliers.
The resulting presentation of \eqref{eq:saddle} takes $A := A_{\cR}$ and $q := q_{\cR}$, with $\cY$ the product of half-lines dual to the rows of $\cR$:
\begin{equation}\label{eq:saddle-split}
    \min_{x\in\cX}\ \max_{y\ge 0}\ \ \langle c, x\rangle \;+\; \langle y,\ q_{\cR}-A_{\cR}x\rangle .
\end{equation}
The standard form is the split $\cE = \emptyset$. 
The same reformulation can also be applied to the dual LP, as in Section~\ref{sec:maxselection}.
We propose to reformulate structured LPs using this absorption split.

\subsection{A useful subset of polyhedra with fast projections}\label{sec:constructs}

Table~\ref{tab:constructs} presents some examples of practical polyhedral structures with efficient Euclidean projections.
\begin{table}[ht]
    \centering
    \small
    \begin{tabular}{@{}l l l l@{}}
        \toprule
        Set & Definition & Projection & Constraint \\
        \midrule
        $\cX$: simplex per group & $\sum_{i\in g}x_i = r_g$, $x \ge 0$ & $O(n)$ & assignment, partition \\
        $\cX$: weighted simplex & $a^\top x = r$, $a > 0$, $x \ge 0$ & $O(n)$ & weighted partition \\
        $\cX$: $\ell_1$ ball & $\norm{x}_1 \le r$ & $O(n)$ & sparsity budget \\
        $\cX$: isotonic cone & $x_1 \le \cdots \le x_n$ & $O(n)$ & monotonicity \\
        $\cY$: simplex per group & $\sum_{j\in G_g}\lambda_{g,j} = c_g$, $\lambda \ge 0$ & $O(n)$ & max-selection \\
        \bottomrule
    \end{tabular}
\caption{Candidate structures for an absorption split, with the constraint family each one absorbs.
The projection column is the best worst-case cost over $n$ coordinates.}
    \label{tab:constructs}
\end{table}

Each row of Table~\ref{tab:constructs} is a family of constraints together with the set it becomes.
Absorption into $\cX$ gives the split presentation \eqref{eq:saddle-split}.

\emph{Assignment and partition rows}, $\sum_{i\in g}x_i = r_g$ with $x \ge 0$ over disjoint groups $g$, give a scaled simplex $\Delta^{r_g}$ per group, projected independently.
Median selection projects a simplex in deterministic linear time \cite{brucker_on_1984,pardalos_algorithm_1990}; sort-and-threshold \cite{held_validation_1974} costs a logarithmic factor more, and expected linear-time pivoting \cite{duchi_efficient_2008} and the variants of \cite{condat_fast_2016} are the fast choices in practice.
These rows back assignment, transportation, covering-by-choice and lot-sizing models.

\emph{A weighted partition}, $a^\top x = r$ with $a > 0$ and $x \ge 0$, is the weighted simplex $\cX = \{x \ge 0 : a^\top x = r\}$, projected in linear time by \cite{perez_efficient_2022}.

\emph{A sparsity budget}, $\norm{x}_1 \le r$, gives $\cX = \{x : \norm{x}_1 \le r\}$, which is a simplex after a sign substitution and inherits its projection.

\emph{Monotonicity}, $x_1 \le \cdots \le x_n$ within a group, gives the isotonic cone, the structure behind staircase models in scheduling and inventory, projected in linear time by pool-adjacent-violators \cite{best_active_1990}.

\emph{Max-selection rows}, $z_g \ge x_j$ over a group $G_g$ whose auxiliary variable $z_g$ carries the group's maximum at optimality, contribute their simplex to $\cY$ instead. Pricing the group with $\lambda_g \ge 0$ makes the free column of $z_g$ an equality $\sum_{j\in G_g}\lambda_{g,j} = c_g$, so the dual block is a scaled simplex per group and the $z$ variables leave the formulation.
Section~\ref{sec:maxselection} carries this out.

\subsection{Max-selection}\label{sec:maxselection}

Max-selection is common in practical LPs: auxiliary variables that select maxima over groups of the $x$ variables at optimality.
Its dual is an assignment family, so it can be absorbed.
Suppose the LP has the form
\begin{equation}\label{eq:maxsel}
    \min_{x \in \cX_0,\ z}\ \langle h, x\rangle + \sum_{g \in G}c_g z_g
    \qquad\text{s.t.}\quad
    \begin{aligned}
     &A_{\cR}x \ge q_{\cR},\\
        &z_g \ge x_j, \ \forall j \in G_g,\ \forall g \in G,
    \end{aligned}
\end{equation}
where each group $G_g \subseteq [n]$ is a set of coordinates, $c_g \ge 0$, and $A_{\cR}x \ge q_{\cR}$ collects whatever rows are left.

\textbf{Reformulate \eqref{eq:maxsel} as a saddle problem.}
\begin{equation}\label{eq:maxsel-saddle}
    \min_{x\in\cX_0}\ \max_{\substack{\lambda_g\in\Delta^{c_g},\ g \in G\\ y_{\cR}\ge0}}\
    \langle h, x\rangle \;+\; \sum_{g\in G}\lambda_g^\top x_{G_g} \;+\; \langle y_{\cR},\ q_{\cR}-A_{\cR}x\rangle ,
\end{equation}
\textbf{which is \eqref{eq:saddle} with $\cX = \cX_0$ and no $z$ variables.}
The rest of this section derives it.

\textbf{The dual of a max-selection family is an assignment family.}
Each primal variable contributes one row to the dual.
Price $z_g \ge x_j$ with $\lambda_{g,j} \ge 0$ and $A_{\cR}x \ge q_{\cR}$ with $y_{\cR} \ge 0$, so that
\[
    L(x,z;\lambda,y_{\cR}) \;=\; \langle h,x\rangle
    \;+\; \sum_{g\in G} z_g\Big(c_g - \sum_{j\in G_g}\lambda_{g,j}\Big)
    \;+\; \sum_{g\in G}\lambda_g^\top x_{G_g}
    \;+\; \langle y_{\cR},\ q_{\cR}-A_{\cR}x\rangle .
\]
Each $z_g$ is free and occurs in no other row, so the row it contributes to the dual is the equality
\[
    \sum_{j\in G_g}\lambda_{g,j} \;=\; c_g ,
    \qquad \lambda_g \ge 0 ,
\]
one per group, over disjoint groups.
That is the assignment family of Table~\ref{tab:constructs}, with simplices of total $c_g$.

\textbf{So absorb it.}
Definition~\ref{def:split} applied to those rows puts each $\lambda_g$ in its own scaled simplex $\Delta^{c_g}$, so the dual projection separates into one simplex projection per group.
Every $z_g$ is then priced out of \eqref{eq:maxsel}, and what is left is \eqref{eq:maxsel-saddle}.

When $\cR = \emptyset$ and $\cX_0$ was itself absorbed, every row of \eqref{eq:maxsel} and of its dual has been absorbed, and we call the split \emph{complete}. Both case studies are complete splits.

\section{Restarted Halpern PDHG and theoretical guarantees}\label{sec:theory}

This section presents the restarted reflected Halpern algorithm the implementations run \cite{lu_restarted_2024}.
Throughout, \eqref{eq:saddle} is a feasible and bounded LP presented under an absorption split, with $\cX \subseteq \RR^n$ and $\cY \subseteq \RR^m$ nonempty polyhedra carrying computable projections, and the diagonal scalings $\Lambda_{\cX},\Lambda_{\cY}$ of \eqref{eq:pdhg} come from the split.
Write $\cZ := \cX\times\cY$, $v = (x,y)$ and $\Lambda := \diag(\Lambda_{\cX},\Lambda_{\cY})$, and let $\cZ^\star$ be the saddle set of $L$ on $\cZ$, which is then nonempty and satisfies $\cZ^\star = \cX^\star\times\cY^\star$ with $\cX^\star,\cY^\star$ the optimal faces of the two value problems.
Write $L_P(x) := \max_{y\in\cY}L(x,y)$ and $L_D(y) := \min_{x\in\cX}L(x,y)$ for the two
value functions, and $\gp(v) := L_P(x) - L_D(y)$ for the duality gap.

\subsection{The algorithm}\label{sec:alg}

One step of PDHG \eqref{eq:pdhg} is a proximal-point step in the canonical norm
\begin{equation}\label{eq:canonical}
    P := \begin{bmatrix}
        \dfrac{\omega}{\eta}\Lambda_{\cX} & A^\top \\[4pt]
        A & \dfrac{1}{\eta\omega}\Lambda_{\cY}
    \end{bmatrix},
\qquad \norm{v}_P := \sqrt{v^\top Pv} ,
\end{equation}
taken on the saddle operator $F(v) := (c - A^\top y + N_{\cX}(x),\ Ax - q + N_{\cY}(y))$, with $N_S$ the normal cone; $F$ is maximally monotone with $F^{-1}(0) = \cZ^\star$.
In this norm $T_\omega$ is the proximal step of $F$, the standard proximal-point
reading of PDHG \cite{he_convergence_2012}.
It is firmly nonexpansive in $\norm{\cdot}_P$, so its reflection $2T_\omega-I$ is
nonexpansive \cite[Proposition~4.4]{bauschke_convex_2017}.
Conjugating by $\Lambda^{-1/2}$ gives
\[
    \Lambda^{-1/2}P\Lambda^{-1/2} \;=\; \begin{bmatrix}
        \dfrac{\omega}{\eta}I & \Lambda_{\cX}^{-1/2}A^\top\Lambda_{\cY}^{-1/2} \\[4pt]
        \Lambda_{\cY}^{-1/2}A\Lambda_{\cX}^{-1/2} & \dfrac{1}{\eta\omega}I
    \end{bmatrix},
\]
so $P\succ0$ exactly when $\eta < 1/\cnorm{A}$, at every $\omega>0$, and
\begin{equation}\label{eq:lammax}
    \norm{v}_P \;\le\; \sqrt{\frac{\omega+\omega^{-1}}{\eta}}\,\norm{v}_\Lambda,
    \qquad\text{and at }\omega = 1,\quad
    \norm{v}_P \;<\; \sqrt{\frac2\eta}\,\norm{v}_\Lambda .
\end{equation}

\paragraph{Halpern iteration.}
Given an anchor $v^{\mathrm a}$, the Halpern iteration \cite{halpern_fixed_1967} applied to $T_\omega$ takes the form
\begin{equation}\label{eq:halpern-plain}
    v^{k+1} = \frac{k+1}{k+2}T_\omega(v^k) + \frac{1}{k+2}v^{\mathrm a} .
\end{equation}
Since $T_\omega$ is nonexpansive in $\norm{\cdot}_P$ and has a fixed point, its iterates
satisfy \cite{lieder_convergence_2021}
\[
    \norm{v^k - T_\omega(v^k)}_P \;\le\; \frac{2\dist_P(v^{\mathrm a},\cZ^\star)}{k+1} .
\]

\paragraph{Reflection.}
The reflection $2T_\omega - I$ has the same fixed points as $T_\omega$.
Applying the Halpern iteration to it gives the reflected Halpern iteration \cite{bauschke_convex_2017}
\begin{equation}\label{eq:halpern}
    v^{k+1} = \frac{k+1}{k+2}\pran{2T_\omega(v^k) - v^k} + \frac{1}{k+2}v^{\mathrm a},
\end{equation}
which halves the constant of \eqref{eq:halpern-plain}
\cite[Corollary~1]{lu_restarted_2024}.

\paragraph{Restarting.}
The restarted reflected Halpern iteration \cite{lu_restarted_2024,lu_cupdlpx_2025} resets the anchor when a restart condition is met.
Let $j$ and $k$ denote the outer and inner iteration counters.
Within outer iteration $j$ it generates
\begin{equation}\label{eq:halpern-restart}
    v^{j,k+1} = \frac{k+1}{k+2}\pran{2T_\omega(v^{j,k}) - v^{j,k}} + \frac{1}{k+2}v^{j,0},
\end{equation}
and when a restart is triggered after inner iteration $k_j$, the next outer iteration starts from the corresponding fixed-point step,
\begin{equation}\label{eq:halpern-anchor}
    v^{j+1,0} := T_\omega(v^{j,k_j}) .
\end{equation}
Write $v^{\mathrm a} := v^{j,0}$ for the current anchor.
Every restart trigger and termination test is evaluated at the plain step $w := T_\omega(v)$, which lies in $\cZ$.
Restarts fire on the triple used by cuPDLPx \cite{lu_cupdlpx_2025}.
\emph{Sufficient decay} fires at
$\norm{v - T_\omega(v)}_P \le \beta_{\mathrm{suff}}\norm{v^{\mathrm a} - T_\omega(v^{\mathrm a})}_P$.
\emph{Necessary decay with stall} fires at
$\norm{v - T_\omega(v)}_P \le \beta_{\mathrm{nec}}\norm{v^{\mathrm a} - T_\omega(v^{\mathrm a})}_P$
when $\norm{v - T_\omega(v)}_P$ is above its value at the previous check.
\emph{Artificial} fires once the epoch has run for a $\beta_{\mathrm{art}}$ fraction of
all steps taken so far.
The anchor \eqref{eq:halpern-anchor} is taken at the updated primal weight, where nonexpansiveness holds its residual at or below that of $v^{j,k_j}$ in the same norm.

The primal weight is updated at restarts only, toward the ratio minimizing $\norm{\cdot}_P$ at the observed displacement, in log space and with a decaying smoothing exponent:
\begin{equation}\label{eq:omega-update}
    \ln\omega \leftarrow (1-\theta)\ln\omega + \theta\ln
    \frac{\norm{y - y^{\mathrm a}}_{\Lambda_{\cY}}}{\norm{x - x^{\mathrm a}}_{\Lambda_{\cX}}},
    \qquad
    \theta \leftarrow \theta_{\mathrm{tgt}} + \theta_{\mathrm{dec}}(\theta - \theta_{\mathrm{tgt}}),
\end{equation}
skipping the update if either norm vanishes.
Early restarts therefore adapt $\omega$ aggressively while later ones, when the iterate is close to optimal and the movement ratios are noisy, perturb it only mildly.

\paragraph{Computational implementation.}
The method is \eqref{eq:halpern-restart} applied to the absorbed saddle problem \eqref{eq:saddle-split}, run at the parameters of Table~\ref{tab:hyperparameters} in Appendix~\ref{app:params} and started at
\[
    v^{0,0} = \pran{\proj_{\cX}(0),\ \proj_{\cY}(0)} .
\]
The instance data is scaled to unit norm, the scale at which
Table~\ref{tab:hyperparameters} is tuned, and unscaled at termination.
Every $\ell$ iterations the run evaluates the relative KKT residual $\kappa$ of
\eqref{eq:lp-def} at $T_\omega(v)$ and stops once $\kappa \le \varepsilon$.
Section~\ref{sec:applications} bounds the fixed-point residual and the distance to
$\cZ^\star$; the implementations terminate on the relative KKT residual.

\subsection{Convergence guarantees of \texorpdfstring{r\(^2\)HPDHG}{r2HPDHG}}\label{sec:pdlp-rates}

The reformulations below are measured against PDHG run on the standard form.
The rates to improve on are the ones that method already has.
Both are from \cite{lu_restarted_2024}, which states them for the standard form.
Appendix~\ref{app:restart} records the proof of Proposition~\ref{prop:luyang}(b) for
$T_\omega$ in the canonical norm \eqref{eq:canonical}.

\begin{mydef}[Sharpness modulus]\label{def:alpha}
$\alpha$ is the largest constant with
$\alpha\,\dist_P(v,\cZ^\star) \le \norm{v - T_\omega(v)}_P$ for every $v$.
\end{mydef}

\begin{prop}[{\cite[Corollary~1, and Theorem~2 with the reflection of Section~6]{lu_restarted_2024}}]\label{prop:luyang}
\leavevmode
\begin{enumerate}[label=(\alph*)]
    \item Run \eqref{eq:halpern-restart} from an anchor $v^{\mathrm a}$ without restarting.
    Then for every $k \ge 0$,
    \[
        \norm{v^k - T_\omega(v^k)}_P \;\le\; \frac{\dist_P(v^{\mathrm a},\cZ^\star)}{k+1} .
    \]
    \item Let $\alpha > 0$ be as in Definition~\ref{def:alpha}, and run the method of
    Section~\ref{sec:alg} with $\omega$ held fixed and sufficient decay as its only
    restart trigger, so that epoch $j$ ends once
    \begin{equation}\label{eq:adaptive-restart}
        \norm{v^{j,k} - T_\omega(v^{j,k})}_P \;\le\; e^{-1}\norm{v^{j,0} - T_\omega(v^{j,0})}_P .
    \end{equation}
    Writing $\tau^j$ for the length of epoch $j$, every epoch satisfies
    $\tau^j \le \lceil e/\alpha\rceil$, and for every $j \ge 0$,
    \[
        \dist_P(v^{j+1,0},\cZ^\star)
        \;\le\; \frac{e^{-j}}{\alpha\pran{\tau^0+1}}\,\dist_P(v^{0,0},\cZ^\star) .
    \]
\end{enumerate}
\end{prop}

The rule \eqref{eq:adaptive-restart} runs without $\alpha$, which enters only the bounds.
It is the sufficient-decay trigger of Section~\ref{sec:alg} at
$\beta_{\mathrm{suff}} = e^{-1}$, the value both implementations use.
PRAQ and PRIA add two enhancements to it: an adaptive $\omega$ at each restart, and the
other two restart triggers of Section~\ref{sec:alg}.

Only $\alpha$ is instance-specific, and it is what a reformulation acts on.

A positive $\alpha$ exists whenever $\cX$ and $\cY$ are polytopes: the gap error bound of
\cite{applegate_faster_2023,lu_infimal_2023,fercoq_quadratic_2023} is then finite on
$\cZ$, and Lemma~\ref{lem:sharp} carries it to the residual at every $v \in \RR^{n+m}$.
The quantifier matters, because the reflected iterates leave $\cZ$.
Those arguments give existence alone, so a value for $\alpha$ is instance-specific work.
Section~\ref{sec:qubo-inst} does it for QUBO.

For the standard form the available constant is a Hoffman constant of the KKT system,
$H_{\mathrm{KKT}}$ \cite{hoffman_approximate_1952}.
Take the baseline at $\Lambda = I$ and $\eta \le 1/(2\cnorm{A})$, the setting
\cite{lu_restarted_2024} analyzes.
Write $R := 2\dist_\Lambda(v^{0,0},\cZ^\star)$ for the distance scale, and let $z^\star$ be
an optimal point nearest $v^{0,0}$.
The restarted iterates stay in the ball about the origin of radius
$\rho := R + 2\norm{z^\star}_2$
\cite[Lemmas~1 and~2(i), with the reflection of Section~6]{lu_restarted_2024}, and on that ball
the error bound of \cite[Lemma~5]{applegate_faster_2023} holds with constant
$H_{\mathrm{KKT}}\sqrt{1+4\rho^2}$.
Proposition~\ref{prop:luyang}(a) and \cite[Theorem~2, with the reflection of
Section~6]{lu_restarted_2024} with that constant put the standard form at
\begin{equation}\label{eq:pdlp-sublinear}
    O\!\pran{\frac{\cnorm{A} R}{\varepsilon}}
\end{equation}
applications of $T_\omega$ for a KKT residual $\varepsilon$, and
\begin{equation}\label{eq:pdlp-linear}
    O\!\pran{\cnorm{A}\,H_{\mathrm{KKT}}\sqrt{1+4\rho^2}\,\log\frac{1}{\varepsilon}}
\end{equation}
applications to contract $\dist_\Lambda(v,\cZ^\star)$ by $\varepsilon$.
These are the two columns the case studies fill in.
A split acts on both, by fixing $\cnorm{A}$ and bounding $\cX$ and $\cY$.

\section{Why reformulate}\label{sec:why}

Absorption changes a solver's performance through two independent channels: the cost of an iteration, and the number of iterations required. In this section we explain how reformulation affects these channels.
\subsection{What changes numerically}\label{sec:whynum}

\paragraph{Cost per iteration.}
An iteration of \eqref{eq:pdhg} in the standard form requires two sparse matrix--vector multiplications, one with $A$ and one with $A^\top$. Under a complete split with identity-like coupling, these products disappear entirely: the interaction terms $A^\top y$ and $Ax$ reduce to elementwise accesses to the other variable. The resulting iteration consists of elementwise arithmetic and two projections, which can be implemented with custom GPU kernels.

The difference is important on GPUs. An SpMV kernel is inherently irregular, whereas an elementwise kernel is uniform. SpMV gathers entries of $x$ through an index array, performs only about two floating-point operations per nonzero, incurs irregular memory access, and can suffer from load imbalance when row lengths vary. Moreover, PDHG requires both $Ax$ and $A^\top y$, which either requires storing the matrix in both orientations or using atomic scatter operations for one of the products. By contrast, the absorbed iteration under a complete split can use contiguous memory access. In some structured cases, a customized kernel can similarly replace SpMV, as in the two examples below.

This speedup applies specifically to the interaction term. The tradeoff is that the additional projection cost introduced by the split can become dominant, as it does for the JRP example in Table~\ref{tab:jrp-cost}.

\paragraph{Iteration count.}
Reducing the number of dual variables lowers the cost of every vector operation, but the more consequential benefit is a reduction in the number of iterations.
The constants of Section~\ref{sec:whytheory} bound the saving; the measured counts realize it.
Against cuPDLPx on the original formulation, PRAQ converges in $1.7$--$4.3$ times fewer iterations and PRIA in $13.6$--$52.3$ times fewer, over the sizes where cuPDLPx converges at all.

\subsection{What changes theoretically}\label{sec:whytheory}

A formulation enters the convergence theory of Section~\ref{sec:pdlp-rates} through a short list of constants.
Absorption acts on each separately, so its effect can be attributed term by term.

\paragraph{A data-free step size.}
The admissible step size satisfies $\eta<1/\cnorm{A}$, where $\cnorm{A}$ is the coupling norm induced by the split (Section~\ref{sec:alg}). In the standard formulation, the step size is therefore determined by the norm of the numerical constraint matrix. Diagonal rescaling can reduce this coupling norm to $O(1)$ \cite{pock_diagonal_2011}, but does so by absorbing the data dependence into the scaling.
Under an absorption split, by contrast, the coupling norm can be bounded purely from the structure of the split, yielding a data-free step size that is fixed for the entire problem family. In particular, $A_{\cR}$ is a row submatrix of the full constraint matrix, so $\norm{A_{\cR}}_2\le\norm{A}_2$. In our case studies, $\cnorm{A}\le1$ for QUBO and $\cnorm{A}=\sqrt{2}$ for JRP.

\paragraph{A tighter termination certificate.}
The absorbed constraints and box constraints are satisfied exactly by construction, so the KKT residual needs to measure only the constraints in $\cR$ and the optimality gap. Each additional absorbed constraint family eliminates a feasibility residual that would otherwise need to be controlled numerically, thereby tightening the termination certificate. Under a complete split, all feasibility residuals vanish, and the KKT residual reduces to the gap term.

\paragraph{An explicit distance scale.}
The sublinear rate of Section~\ref{sec:pdlp-rates} is measured against the distance
scale $R$.
In the standard formulation, $\cY$ is a product of half-lines, so $\cZ$ has infinite diameter and $R$ can only be controlled through the distance to the saddle set, which is unknown before solving the problem. 
Under a complete split $\cZ$ is likely a bounded polytope and $R \le 2\diam_\Lambda(\cZ)$ can be fully
explicit.
In our case studies, $\diam_\Lambda(\cZ)^2 \le 2m+2$ for QUBO~\eqref{eq:qubo-constants}
and $\diam_\Lambda(\cZ)^2 \le 2NT+2T$ for the JRP~\eqref{eq:jrp-constants}.

\paragraph{Sharpness constant.}
Restarted PDHG converges linearly at a rate set by a sharpness constant, and every linear program supplies one in every formulation (Section~\ref{sec:pdlp-rates}).
The standard-form analysis needs the Hoffman constant of the whole KKT system.
Absorption removes the multipliers of $\cE$ along with its rows, so the absorbed formulation carries the sharpness modulus of Definition~\ref{def:alpha} for a saddle problem in fewer variables.
Making that modulus \emph{explicit} is a separate question, and for QUBO we bound it in Proposition~\ref{prop:qubo-sharp}.

Table~\ref{tab:constants} summarizes the four constants.

\begin{table}[ht]
    \centering
\small
\renewcommand{\arraystretch}{1.4}
\begin{tabular}{@{}>{\raggedright\arraybackslash}p{2.4cm}>{\raggedright\arraybackslash}p{5.0cm}>{\raggedright\arraybackslash}p{6.3cm}@{}}
\toprule
\textbf{Constant} & \textbf{Standard form} & \textbf{After absorption} \\
\midrule
Coupling $\cnorm{A}$
& The full constraint matrix, unbounded in general, and growing with the largest
  number of rows that share a column.
& Only the unabsorbed rows, and at most $1$ for QUBO, exactly $\sqrt2$ for the JRP. \\
Termination
& KKT residual: primal infeasibility, dual infeasibility, gap.
& One term drops per absorbed family. For QUBO and JRP all feasibility residuals vanish, and the KKT residual reduces to the gap term. \\
Distance scale $R$
& Bounded only through the saddle set, $\cZ$ having infinite diameter.
& $R \le 2\diam_\Lambda(\cZ)$, read off the absorbed sets, once the split is complete. \\
Sharpness
& Global Hoffman constant of the KKT system.
& The sharpness modulus of the reformulated saddle problem: explicit when the reformulated polytope has a lattice
  of vertex values (Appendix~\ref{app:sharp}). \\
\bottomrule
\end{tabular}
\caption{The constants the convergence theory reads off a formulation, before and after a split.}
    \label{tab:constants}
\end{table}

\section{Applications}\label{sec:applications}

In this section we apply the reformulation to two case studies.
Both splits are complete, so $\cX$ and $\cY$ are polytopes, the absorbed rows hold exactly at every termination check, and the certificate is the duality gap alone.
Both are reformulated by the recipe of Section~\ref{sec:maxselection} and instantiated against Section~\ref{sec:theory}, and they differ in whether the sharpness modulus of Definition~\ref{def:alpha} can be made explicit.
Tables~\ref{tab:qubo-rates} and~\ref{tab:jrp-rates} record the theoretical convergence rate improvements against \eqref{eq:pdlp-sublinear} and \eqref{eq:pdlp-linear}. Tables~\ref{tab:qubo-solve} and~\ref{tab:jrp-solve} give computational results.

\subsection{Quadratic unconstrained binary optimization (QUBO)}\label{sec:qubo}
Given a symmetric $W \in \RR^{n\times n}$, quadratic unconstrained binary optimization is
\begin{equation}\label{eq:qubo}
    \min_{x\in\{0,1\}^n}\ x^\top Wx .
\end{equation}
It is NP-hard and is the standard target format for Ising machines and quantum annealers.
Its LP relaxation, the \emph{roof dual} \cite{hammer_roof_1984}, is the workhorse bound inside branch-and-bound and the source of variable fixings.
This section reformulates that relaxation and instantiates Section~\ref{sec:theory} on it, with explicit constants.

\subsubsection{The standard linearization}\label{sec:qubo-lin}

Using $x_i^2 = x_i$ to route the diagonal through the linear term, write
\begin{equation}\label{eq:qubo-split}
    x^\top Wx \;=\; \sum_i d_i x_i + \sum_{i<j} 2W_{ij}x_ix_j,
    \qquad d_i := W_{ii},
\end{equation}
and set $\Ep := \{i<j : W_{ij} > 0\}$, $\Em := \{i<j : W_{ij} < 0\}$, $E := \Ep\cup\Em$ with $m := |E|$, and $w_e := 2|W_{ij}|$.
We read $E$ as the edge set of a graph on $[n]$: one node per variable, one edge per off-diagonal nonzero of $W$, that is per quadratic term of \eqref{eq:qubo-split}.
Write $\deg(i)$ for the number of edges at node $i$.
Introducing $z_e$ for each product and orienting each inequality against the sign of its objective coefficient gives the classical linearization of Glover and Woolsey \cite{glover_converting_1974}, whose relaxation we write $\mathrm{(P)}$:
\begin{mini}{x,\,z}{
    d^\top x + \sum_{e\in\Ep}w_ez_e - \sum_{e\in\Em}w_ez_e
  }{\label{eq:qubo-lp}}{}
  \addConstraint{z_e}{\ge x_i+x_j-1,\quad z_e \ge 0,\quad}{e = \{i,j\}\in\Ep}
  \addConstraint{z_e}{\le x_i,\quad z_e \le x_j,\quad}{e = \{i,j\}\in\Em}
  \addConstraint{x}{\in[0,1]^n,\quad z \in [0,1]^E.}{}
\end{mini}

\paragraph{Observations on the linearization.}
The following observations clarify the role of the variable bounds and diagonal terms in~\eqref{eq:qubo-lp}.
\begin{enumerate}[label=(\roman*)]
    \item \textbf{Lower bounds on $z_e$.}
    For $e\in\Ep$, the positive coefficient $w_e$ makes $z_e = \max\{0,x_i+x_j-1\}$ at optimality.
    The lower bound $z_e\ge0$ is essential: without it, $x_i=x_j=0$ would permit $z_e=-1$, contributing $-w_e$ to the objective instead of $0$.
    For $e\in\Em$, the negative coefficient makes $z_e=\min\{x_i,x_j\}\in[0,1]$ at optimality.
    \item \textbf{Upper bounds on $x$.}
    The bounds $x\le1$ must be retained when relaxing $x\in\{0,1\}^n$.
    For $e=\{i,j\}\in\Ep$, combining $z_e\ge x_i+x_j-1$ with $z_e\le1$ gives $x_i+x_j\le2$, so each variable needs its own bound.
    \item \textbf{Treatment of diagonal terms.}
    The binary identity $x_i^2=x_i$ gives the linear contribution $W_{ii}x_i$ in~\eqref{eq:qubo-split}.
    Introducing a separate product variable $z_{ii}$ can weaken the relaxation: when $W_{ii}>0$, this replaces $W_{ii}x_i$ by $W_{ii}\max\{0,2x_i-1\}$ at optimality.
    At $x_i=\tfrac12$, the latter is $0$, whereas the linear contribution is $W_{ii}/2$.
\end{enumerate}

\subsubsection{The reformulation}\label{sec:qubo-reform}
\paragraph{The reformulation.} QUBO has the max-selection structure of Section~\ref{sec:maxselection} on both edge families, with groups of size two.
On $\Ep$ the coefficient is $+w_e \ge 0$ and at optimality $z_e = \max\{0, x_i+x_j-1\}$, a weighted maximum over the group $(0,\ x_i+x_j-1)$.
On $\Em$ the coefficient is $-w_e \le 0$ and at optimality $z_e = \min\{x_i,x_j\}$. Absorbing the sign into the maximum, $-w_e\min\{x_i,x_j\} = w_e\max\{-x_i,-x_j\}$, puts $\Em$ in the same nonnegatively weighted form, with group $(-x_i,\ -x_j)$ and weight $w_e$.
Eliminating $z$ therefore leaves $\mathrm{(P)} = \min_{x\in[0,1]^n}F(x)$ for the convex piecewise linear
\begin{equation}\label{eq:qubo-F}
    F(x) \;:=\; d^\top x + \sum_{e = \{i,j\} \in\Ep}w_e\max\{0,\ x_i+x_j-1\}
    \;+\; \sum_{e = \{i,j\} \in\Em}w_e\max\{-x_i,\ -x_j\} .
\end{equation}
Write $c_k := d_k - \sum_{e\in\Em:\,k=\max e} w_e$, let $q_e := -1$ on $\Ep$ and $0$ on $\Em$, let
\begin{equation}\label{eq:qubo-Y}
    \cY \;:=\; \prod_{e\in E}\,[0,\ w_e],
\end{equation}
and let $\tilde A : \RR^n\to\RR^E$ be the signed incidence operator
\begin{equation}\label{eq:qubo-incidence}
    (\tilde Ax)_e \;:=\;
    \begin{cases}
        x_i - x_j, & e = \{i,j\}\in\Em,\\
        -(x_i + x_j), & e = \{i,j\}\in\Ep,
    \end{cases}
    \qquad i<j .
\end{equation}

\begin{prop}[The QUBO reformulation]\label{prop:qubo-reform}
$\mathrm{(P)}$ is the saddle problem
\begin{equation}\label{eq:qubo-bilinear}
    \min_{x\in[0,1]^n}\ \max_{\lambda\in\cY}\ L(x,\lambda)
    \;:=\; \langle c,\ x\rangle + \langle \lambda,\ q - \tilde Ax\rangle ,
\end{equation}
which is \eqref{eq:saddle} with $\cX = [0,1]^n$: the inner maximum is $F(x)$ at every $x \in [0,1]^n$, so the two have the same optimal value, attained at the same $x^\star$.
\end{prop}

\begin{proof}
We prove Proposition~\ref{prop:qubo-reform} in three steps: price the two rows of each edge, absorb the simplex that pricing produces, and evaluate the inner maximum.

\emph{Step 1: the dual of each edge is a simplex.}
Write $v_e$ for the variable carrying the maximum of edge $e$ --- $z_e$ on $\Ep$ and $-z_e$ on $\Em$ --- so that in both cases $\mathrm{(P)}$ pays $w_e v_e$ against the two rows of the group.
Pricing those rows with $M_{e,1}, M_{e,2} \ge 0$ leaves
\[
    e\in\Em:\ \ v_e\pran{w_e - M_{e,1} - M_{e,2}} - M_{e,1}x_i - M_{e,2}x_j,
    \qquad
    e\in\Ep:\ \ v_e\pran{w_e - M_{e,1} - M_{e,2}} + M_{e,2}\pran{x_i+x_j-1},
\]
and each $v_e$ is free and occurs in no other row, so by Section~\ref{sec:maxselection} the row its column contributes to the dual is
\begin{equation}\label{eq:qubo-dualrow}
    M_{e,1} + M_{e,2} \;=\; w_e, \qquad M_e \ge 0 ,
\end{equation}
one per edge, over disjoint groups. This is the assignment family of Table~\ref{tab:constructs}, with a simplex of total $w_e$.

\emph{Step 2: absorb the simplices.}
Definition~\ref{def:split} applied to those rows puts each $M_e$ in its own scaled simplex $\Delta^{w_e}$, one per edge.
On $\Ep$ the group member $0$ carries no $x$, so $M_{e,1}$ reaches the objective only through \eqref{eq:qubo-dualrow}; write $\lambda_e := M_{e,2}$ for the multiplier that survives.
What is left of $\mathrm{(P)}$ is a saddle problem in the original variables:
\begin{equation}\label{eq:qubo-saddle-x}
    \min_{x\in[0,1]^n}\ \max_{M_e\in\Delta^{w_e}}\quad
    \langle d,\ x\rangle
    \;-\!\!\sum_{e=\{i,j\}\in\Em}\!\!\pran{M_{e,1}x_i + M_{e,2}x_j}
    \;+\!\!\sum_{e=\{i,j\}\in\Ep}\!\!\lambda_e\pran{x_i + x_j - 1}.
\end{equation}
Each $\Delta^{w_e}$ is two-dimensional, hence a segment with a single degree of freedom, so a change of variables collapses the two dual coordinates of an edge to one,
\begin{equation}\label{eq:qubo-collapse}
    e \in \Em:\ \ M_e = (\lambda_e,\ w_e-\lambda_e),
    \qquad
    e \in \Ep:\ \ M_e = (\lambda_e,\ \lambda_e),
    \qquad \lambda_e \in [0,\ w_e],
\end{equation}
and the dual set becomes the box \eqref{eq:qubo-Y}.
Substituting \eqref{eq:qubo-collapse} into \eqref{eq:qubo-saddle-x} leaves the same saddle problem with one dual coordinate per edge:
\begin{equation}\label{eq:qubo-saddle-lambda}
    \min_{x\in[0,1]^n}\ \max_{\lambda\in\cY}\quad
    \langle d,\ x\rangle
    \;-\!\!\sum_{e=\{i,j\}\in\Em}\!\!\pran{\lambda_e x_i + (w_e-\lambda_e)x_j}
    \;+\!\!\sum_{e=\{i,j\}\in\Ep}\!\!\lambda_e\pran{x_i + x_j - 1}.
\end{equation}
On $\Em$ the term $-w_ex_j$ carries no dual variable, so it belongs to the primal cost, which is what $c$ collects; gathering the rest against $q$ and $\tilde A$ gives \eqref{eq:qubo-bilinear}.

\emph{Step 3: the inner maximum is $F$.}
Substitute \eqref{eq:qubo-collapse} into \eqref{eq:qubo-saddle-x}; everything but $\langle d,x\rangle$ separates over edges.
On $e=\{i,j\}\in\Em$ the contribution is $-\lambda_ex_i - (w_e-\lambda_e)x_j = -w_ex_j + \lambda_e(x_j-x_i)$: the part carrying no dual variable is absorbed into $c$, and the rest is $\lambda_e(q_e - (\tilde Ax)_e)$.
Its maximum over $[0,w_e]$ is $-w_ex_j + w_e\max\{0,x_j-x_i\}$, which is $-w_ex_i$ when $x_i \le x_j$ and $-w_ex_j$ otherwise, i.e.\ $w_e\max\{-x_i,-x_j\}$ in both cases.
On $e\in\Ep$ the contribution is $\lambda_e(x_i+x_j-1) = \lambda_e(q_e - (\tilde Ax)_e)$, whose maximum over $[0,w_e]$ is $w_e\max\{0,x_i+x_j-1\}$.
Summing gives $\max_{\lambda\in\cY}L(x,\lambda) = F(x)$ by \eqref{eq:qubo-F}, and minimizing both sides over $[0,1]^n$ gives the claim.
\end{proof}

\paragraph{Diagonal scaling.}
$\Lambda_{\cX}$ weights node $i$ by $\deg(i)$, and
$\Lambda_{\cY}$ weights every edge by $2$.
These are the diagonal preconditioner of Pock and Chambolle \cite{pock_diagonal_2011} at
exponent $1$, since $\sum_e|\tilde A_{ei}|=\deg(i)$ and $\sum_i|\tilde A_{ei}|=2$.
In~\eqref{eq:pdhg} the primal step at node $i$ is scaled by $1/\deg(i)$ and the dual step
by $1/2$.
Both box projections remain coordinatewise under these positive diagonal scalings.

\paragraph{Step size.}
Under these scalings the coupling operator satisfies
\begin{equation}\label{eq:qubo-Anorm}
    \cnorm{\tilde A} \;=\; \norm{\Lambda_{\cY}^{-1/2}\tilde A\Lambda_{\cX}^{-1/2}}_2
    \;\le\; 1.
\end{equation}
This is the Pock--Chambolle bound at exponent $1$ \cite{pock_diagonal_2011}.
Any $0<\eta<1$ is therefore admissible at every primal weight, independently of the graph
size and maximum degree.
We use $\eta=0.999$ in the implementation.

\paragraph{Per-iteration cost.}
With the scalings specified above, the PDHG updates~\eqref{eq:pdhg} take the form
\begin{equation}\label{eq:qubo-pdhg}
\begin{aligned}
    x^{k+1} &:= \proj_{[0,1]^n}\!\pran{x^k - \tfrac{\eta}{\omega}\Lambda_{\cX}^{-1}\pran{c - \tilde A^\top\lambda^k}},\\
    \lambda^{k+1} &:= \proj_{\cY}\!\pran{\lambda^k + \tfrac{\eta\omega}{2}\pran{q - \tilde A(2x^{k+1}-x^k)}},
\end{aligned}
\end{equation}
Each iteration requires two box projections and one application each of $\tilde A$ and $\tilde A^\top$.
The projections are evaluated coordinatewise, while the two operator applications are implemented as gather and scatter operations over the edge list.
The resulting arithmetic cost is $O(n+m)$ per iteration.
The nonzero entries of $\tilde A$ are $\pm1$ and fixed by the signs of $W_{ij}$, so the edge list carries the whole coupling.
PDHG applied directly to~\eqref{eq:qubo-lp} takes two sparse matrix--vector products per iteration against the stored constraint matrix.

\paragraph{Termination.}
The QUBO reformulation allows termination to be checked using the primal--dual gap alone.
At each scheduled check, we evaluate the projected PDHG point $(x,\lambda):=T_\omega(v)$, which satisfies $x\in[0,1]^n$ and $\lambda\in\cY$, as described in Section~\ref{sec:alg}.
The value functions $L_P$ and $L_D$ of Section~\ref{sec:theory} bound the optimal value
from above and below.
For QUBO both are available in closed form. Writing $g(\lambda):=c-\tilde A^\top\lambda$, we obtain
\begin{equation}\label{eq:qubo-bounds}
    L_D(\lambda) = \sum_{i=1}^n\min\{0,\ g_i(\lambda)\} - \sum_{e\in\Ep}\lambda_e,
    \qquad
    L_P(x) = F(x).
\end{equation}
The lower bound follows by minimizing $L(x,\lambda)=\langle g(\lambda),x\rangle+\langle\lambda,q\rangle$ over the primal box; the upper bound follows from Proposition~\ref{prop:qubo-reform} by recovering the eliminated product variables.

We terminate when the relative duality gap satisfies
\[
    \frac{L_P(x)-L_D(\lambda)}{1+|L_P(x)|+|L_D(\lambda)|}\le\varepsilon,
\]
where the objective values in this test are evaluated at the original data scale.
The monitored point satisfies both box constraints and the eliminated variables are recovered exactly, so the gap alone certifies.

\subsubsection{Theoretical improvements}\label{sec:qubo-inst}

\begin{ass}[Integral data]\label{ass:qubo-integral}
$d \in \ZZ^n$ and $w \in \ZZ^E$.
\end{ass}

\begin{thm}[Convergence of PRAQ]\label{thm:qubo-rate}
For the roof dual \eqref{eq:qubo-lp} of an instance $W$ satisfying
Assumption~\ref{ass:qubo-integral}, with every node carrying at least one edge
and $n \ge 5$, run the method of Section~\ref{sec:alg} on
\eqref{eq:qubo-bilinear} at a fixed $\omega = 1$ and $0 < \eta < 1$, on the
normalized instance $(d;w)/\norm{(d;w)}_2$.
Write $v^k$ for its iterates without restarts, $v^{j,k}$ for those of its $j$th
epoch under restarts, and $\norm{\cdot}_P$ for the canonical norm
\eqref{eq:canonical} at that primal weight.
\begin{enumerate}[label=(\alph*)]
    \item \textbf{Sublinear convergence without restarts.}
    For every $k \ge 0$,
    \[
        \norm{v^k - T_\omega(v^k)}_P \;\le\; \frac{2}{k+1}\sqrt{\frac{2m}{\eta}} .
    \]
    \item \textbf{Linear convergence with adaptive restarts.}
    Restarting by sufficient decay \eqref{eq:adaptive-restart}, every epoch has length at most
    $\left\lceil e\pran{1 + \tfrac{16}{\eta}\norm{W}_F\,m\sqrt n}\right\rceil$, and for
    every $j \ge 0$,
    \[
        \dist_P(v^{j+1,0},\cZ^\star) \;\le\;
        e^{-j}\pran{1 + \frac{16}{\eta}\norm{W}_F\,m\sqrt n}\sqrt{\frac{2m}{\eta}} .
    \]
\end{enumerate}
\end{thm}

\begin{proof}
The reformulation has
\begin{equation}\label{eq:qubo-constants}
    \cnorm{\tilde A} \;\le\; 1,
    \qquad
    \diam_\Lambda(\cZ)^2 \;=\; 2m + \frac{2\norm{w}_2^2}{\norm{(d;w)}_2^2} \;\le\; 2m+2,
\end{equation}
by \eqref{eq:qubo-Anorm}, and by $\Lambda_{\cX}$ weighting node $i$ by $\deg(i)$ and
$\Lambda_{\cY}$ weighting each edge by $2$ with $\sum_i\deg(i) = 2m$ and
$\norm{w}_2 \le \norm{(d;w)}_2$.
The sharpness modulus satisfies
\begin{equation}\label{eq:qubo-alpha}
    \frac1\alpha \;\le\; 1 + \frac{16}{\eta}\norm{W}_F\,m\sqrt n
\end{equation}
by Proposition~\ref{prop:qubo-sharp}, proved in Appendix~\ref{app:sharp}.
The initial point lies in $\cZ$, so \eqref{eq:lammax} at $\omega = 1$ and $m \ge 1$ give
\[
    \dist_P(v^{0,0},\cZ^\star) \;<\; \sqrt{\frac{2}{\eta}}\,\diam_\Lambda(\cZ)
    \;\le\; 2\sqrt{\frac{2m}{\eta}} .
\]
Part (a) is Proposition~\ref{prop:luyang}(a) with this bound, part (b)
Proposition~\ref{prop:luyang}(b) with this $\alpha$ and $\tau^0 \ge 1$.
\end{proof}

\begin{table}[ht]
    \centering
    \small
    \begin{tabular}{@{}l l l l@{}}
        \toprule
        & Sublinear rate & Linear rate & Iteration cost \\
        \midrule
        PDLP formulation & $O\!\pran{\cnorm{A}R/\varepsilon}$
          & $O\!\pran{\cnorm{A}H_{\mathrm{KKT}}\sqrt{1+4\rho^2}\log(1/\varepsilon)}$ & $O(n+m)$ \\
        Reformulation & $O\!\pran{\sqrt m/\varepsilon}$
          & $O\!\pran{\norm{W}_F\,m\sqrt{n}\,\log(1/\varepsilon)}$ & $O(n+m)$ \\
        \bottomrule
    \end{tabular}
\caption{Applications of $T_\omega$ to reach accuracy $\varepsilon$, over $n$ nodes and
$m$ edges, with arithmetic per iteration in the last column.
The first row is \eqref{eq:pdlp-sublinear} and \eqref{eq:pdlp-linear}, in the KKT
residual and in the contraction of $\dist_\Lambda(v,\cZ^\star)$.
The second is Theorem~\ref{thm:qubo-rate}, in the $\norm{\cdot}_P$ fixed-point residual
and in $\dist_P(v,\cZ^\star)$.}
    \label{tab:qubo-rates}
\end{table}
Table~\ref{tab:qubo-rates} presents the effect of reformulation on the constants governing convergence.
For the standard PDLP formulation, the sublinear bound involves the full coupling matrix through $\cnorm{A}$ and the initial distance $R$. While they are potentially bounded by the dimension of the problem and information about the initial solution level-set, the bound is loose.
For the reformulation, the diagonal scaling gives $\cnorm{\tilde A}\le1$, and the bounded primal and dual boxes yield an explicit bound on the initial distance in terms of $m$, providing an explicit and tight sublinear rate. Interestingly, this sublinear bound is independent from the dimension of the problem $n$, only on the number of edges $m$. Furthermore, we can also bound the linear convergence rate. For original PDLP formulation, the linear rate bound additionally depends on $H_{\mathrm{KKT}}$ and the radius $\rho$, which are data dependent and hard to characterize.
Under the assumptions of Theorem~\ref{thm:qubo-rate}, half-integrality of both value polytopes further yields an explicit sharpness estimate.
The problem-dependent terms in the resulting rates are therefore expressed using $m$, $n$, and $\norm{W}_F$, without an unquantified global Hoffman constant.
The different residuals and norms specified in the caption must be taken into account when comparing these bounds.

Reformulation also changes the operations within each iteration.
Although both formulations have $O(n+m)$ arithmetic cost, the reformulation eliminates the product variables and evaluates the signed incidence operator directly through gather and scatter operations over the edge list.
Together with coordinatewise box projections, this structure permits a GPU implementation that avoids storing and multiplying by the full LP constraint matrix.
It can reduce matrix storage and memory traffic even when the asymptotic per-iteration cost is unchanged.

\subsubsection{Computational results}\label{sec:qubo-results}

\paragraph{Implementation.}
The product $\tilde Ax$ is computed by gathering the primal values at the two endpoints of each edge and forming their signed combination.
The transpose product $\tilde A^\top\lambda$ is computed by scattering signed dual contributions to the edge endpoints and accumulating them at each node.
The primal and dual box projections use branch-free minimum and maximum operations to clip each coordinate to $[0,1]$ and $[0,w_e]$, respectively.
The remaining iterate updates are elementwise vector operations.

\paragraph{Instances.}
Three generated families, $2282$ instances with varying edge counts.
Given $n$, a random instance draws $10n$ uniformly random index pairs and keeps the distinct off-diagonal ones, an average degree of $20$.
Edge weights are log-normal with median $1$ and log-standard deviation $2$, negated with probability $\tfrac14$.
Each diagonal entry is a tenth of the node's incident absolute weight times a draw from $\mathrm{Uniform}(-1,1)$.

\paragraph{Protocol.}
Every solver that takes a tolerance is asked for a relative duality gap of $10^{-8}$ and a relative feasibility tolerance of $10^{-9}$, under a five-minute limit per instance.
PRAQ and cuPDLPx \cite{lu_cupdlpx_2025} run on the same NVIDIA H100 80\,GB, one instance at a time.
The third baseline computes the roof dual as a maximum flow with OR-Tools' \textsf{SimpleMaxFlow} \cite{perron_ortools_2024} on CPU.
Those three each solved all $2282$ instances.
Gurobi \cite{gurobi_optimization_llc_gurobi_2026} was run on the same instances on $16$ CPU cores under the same limit, solving $1954$ instances.

\paragraph{Results.}
Table~\ref{tab:qubo-solve} is the running time comparison.
PRAQ gains on cuPDLPx and Gurobi in solve time at every instance size, and is faster than the maximum flow solver on all but the smallest instances.

\begin{table}[ht]
    \centering
    \small
    \begin{tabular}{@{}l r rrr rr r@{}}
        \toprule
        & & & & & \multicolumn{2}{c}{Gurobi} & \\
        \cmidrule(lr){6-7}
        Edges & Instances & PRAQ & cuPDLPx & Max-flow & s & solved & cuPDLPx/PRAQ \\
        \midrule
        $10^3$--$10^4$   & 897 & 0.0054 & 0.052 & \textbf{0.0026} & 0.10  & 897 & $9.6$ \\
        $10^4$--$10^5$   & 767 & \textbf{0.012} & 0.080 & 0.024 & 1.4   & 767 & $6.8$ \\
        $10^5$--$10^6$   & 385 & \textbf{0.048} & 0.36 & 0.47 & $83$ & 290 & $7.5$ \\
        $10^6$--$10^7$   & 193 & \textbf{0.59} & 3.7 & 13.2 & ---   & 0   & $6.4$ \\
        $10^7$--$10^8$   & 40  & \textbf{3.7} & 28.4 & 96.6 & ---   & 0   & $7.7$ \\
        \bottomrule
    \end{tabular}
\caption{Solve time in seconds, shifted (by $1$\, ms) geometric mean by edge count.
The paired column gives how many instances Gurobi finished within the limit; an unfinished bucket's entry is a floor. PRAQ, cuPDLPx and max-flow finished all instances.}
    \label{tab:qubo-solve}
\end{table}

Tables~\ref{tab:qubo-cost} and~\ref{tab:qubo-iters} show that the running-time advantage of PRAQ comes from both cheaper iterations and fewer iterations.
Across the five size groups, PRAQ reduces the median time per iteration by factors of $2.3$--$4.8$ and the iteration count by factors of $1.7$--$4.3$ relative to cuPDLPx.
The lower iteration cost is consistent with evaluating the coupling through custom gather and scatter operations over the edge list, together with coordinatewise box projections.
The reduction in iteration counts provides an additional benefit beyond the lower computational cost of each step, illustrating the two mechanisms discussed in Section~\ref{sec:whynum}.

Table~\ref{tab:qubo-memory} shows a corresponding reduction in peak device memory.
For the two largest size groups, cuPDLPx uses $3.0$--$3.5$ times as much device memory as PRAQ.
In particular, for instances with $10^7$--$10^8$ edges, PRAQ uses $2.26$\,GiB compared with $7.83$\,GiB for cuPDLPx.
This is consistent with eliminating the auxiliary product variables and representing the coupling by an edge list instead of storing the full LP constraint matrix.

\begin{table}[ht!]
    \centering
    \small
    \begin{minipage}[t]{0.48\linewidth}
    \centering
    \begin{tabular}{@{}l rr r@{}}
        \toprule
        Edges & PRAQ & cuPDLPx & cuPDLPx/PRAQ \\
        \midrule
        $10^3$--$10^4$ & \textbf{13.8} & 31.6 & $2.3$ \\
        $10^4$--$10^5$ & \textbf{13.9} & 39.1 & $2.8$ \\
        $10^5$--$10^6$ & \textbf{36.3} & 173 & $4.8$ \\
        $10^6$--$10^7$ & \textbf{353} & 1130 & $3.2$ \\
        $10^7$--$10^8$ & \textbf{1760} & 4400 & $2.5$ \\
        \bottomrule
    \end{tabular}
\caption{Median microseconds per iteration.}
    \label{tab:qubo-cost}
    \end{minipage}\hfill
    \begin{minipage}[t]{0.48\linewidth}
    \centering
    \begin{tabular}{@{}l rr r@{}}
        \toprule
        Edges & PRAQ & cuPDLPx & cuPDLPx/PRAQ \\
        \midrule
        $10^3$--$10^4$ & \textbf{376} & 1617 & $4.3$ \\
        $10^4$--$10^5$ & \textbf{833} & 1952 & $2.3$ \\
        $10^5$--$10^6$ & \textbf{1334} & 2309 & $1.7$ \\
        $10^6$--$10^7$ & \textbf{1720} & 3471 & $2.0$ \\
        $10^7$--$10^8$ & \textbf{2023} & 6385 & $3.2$ \\
        \bottomrule
    \end{tabular}
\caption{Iterations to converge, shifted geometric mean with shift $10$.}
    \label{tab:qubo-iters}
    \end{minipage}
\end{table}

\begin{table}[ht!]
    \centering
    \small
    \begin{tabular}{@{}l r rr r@{}}
        \toprule
        Edges & Instances & PRAQ & cuPDLPx & cuPDLPx/PRAQ \\
        \midrule
        $10^6$--$10^7$ & 193 & \textbf{1.39\,GiB} & 4.22\,GiB & $3.0$ \\
        $10^7$--$10^8$ & 40  & \textbf{2.26\,GiB} & 7.83\,GiB & $3.5$ \\
        \bottomrule
    \end{tabular}
\caption{Largest device memory held during a solve, sampled with NVML.}
    \label{tab:qubo-memory}
\end{table}

\subsection{The joint replenishment problem (JRP)}\label{sec:jrp}

The joint replenishment problem (JRP) is a classical model in deterministic inventory theory, studied since the 1960s as the multi-item extension of the Wagner--Whitin lot-sizing problem \cite{wagner_dynamic_1958}. A review of its literature and industrial applications is in \cite{khouja_review_2008}.
Its stochastic variant remains an active computational target in high dimensions \cite{ata_computational_2026}.
Over a horizon of $T$ periods, each of $N$ items must have its per-period demand served by orders placed in earlier periods.
Ordering item $i$ in a period costs $K_i$. 
Any period in which \emph{some} item is ordered additionally costs $K_0$, modelling a shared fixed cost such as a truck dispatch.
Serving the demand of item $i$ in period $t$ from an order in period $s\le t$ costs $H_{ist}$.
All costs are nonnegative.
The joint cost couples the items and makes the problem hard. 
The JRP is NP-hard \cite{arkin_computational_1989} and remains APX-hard in natural special cases \cite{nonner_approximating_2009}.

The strongest approximation algorithms for the JRP are LP-based.
Levi, Roundy and Shmoys gave a primal--dual $2$-approximation and showed that the natural LP relaxation has integrality gap at most $2$ \cite{levi_primal-dual_2004,levi_primal-dual_2006}.
LP-rounding refinements improved the ratio to $1.8$ \cite{levi_improved_2006} and then $1.791$ \cite{bienkowski_better_2014}.
Every one of these algorithms must first solve the relaxation, so the relaxation is the computational bottleneck at scale: it has $\Theta(NT^2)$ variables and constraints, reaching over a hundred million nonzeros for a few hundred items and periods.

\subsubsection{The LP relaxation}\label{sec:jrp-lp}

The variable $x_{ist}$ is the fraction of item $i$'s period-$t$ demand served from an order in period $s\le t$, $y_{is}$ indicates an order for item $i$ in period $s$, and $z_s$ a joint order in period $s$:
\begin{mini}{x,\,y,\,z}{
    \sum_{s} K_0z_s + \sum_{i,s}K_iy_{is} + \sum_{i,t,s\le t}H_{ist}x_{ist}
  }{\label{eq:jrp-lp}}{}
  \addConstraint{\textstyle\sum_{s\le t}x_{ist}}{= 1,\quad}{i\in[N],\ t\in[T]}
  \addConstraint{x_{ist}-y_{is}}{\le 0,\quad}{i\in[N],\ t\in[T],\ s\le t}
  \addConstraint{x_{ist}-z_s}{\le 0,\quad}{i\in[N],\ t\in[T],\ s\le t}
  \addConstraint{x_{ist},\ y_{is},\ z_s}{\ge 0.}{}
\end{mini}
The assignment rows force every demand to be served; the covering rows charge the item and joint order costs whenever an assignment uses period $s$.
Because $K, K_0 \ge 0$, at optimality
\[
y_{is} = \max_{t\ge s}x_{ist}, \qquad z_s = \max_{i,\ t\ge s}x_{ist},
\]
which is the max-selection structure of Section~\ref{sec:maxselection} with two dual blocks, whose groups are indexed by $(i,s)$ and by $s$.
The method is run on $(\vc(H);K;K_0)/\norm{(\vc(H);K;K_0)}_2$.

\subsubsection{The reformulation}\label{sec:jrp-reform}

Writing $\RR_u^{N\times T\times T}$ for the nonnegative tensors supported on the upper triangle of their last two indices, let
\begin{align*}
    \Delta^1_{1,3} &:= \left\{X \in \RR_u^{N\times T\times T} \;\middle|\;
      \textstyle\sum_{s=1}^{t}X_{ist} = 1,\ i\in[N],\ t\in[T]\right\}, \\
    \Delta^K_{1,2} &:= \left\{\tY \in \RR_u^{N\times T\times T} \;\middle|\;
      \textstyle\sum_{t=s}^{T}\tY_{ist} = K_i,\ i\in[N],\ s\in[T]\right\}, \\
    \Delta^{K_0}_{2} &:= \left\{\tZ \in \RR_u^{N\times T\times T} \;\middle|\;
      \textstyle\sum_{i=1}^{N}\sum_{t=s}^{T}\tZ_{ist} = K_0,\ s\in[T]\right\},
\end{align*}
the superscript giving each simplex's total and the subscript the indices held fixed
in the summation.

\begin{prop}[The JRP reformulation]\label{prop:jrp-reform}
Problem \eqref{eq:jrp-lp} is the saddle problem
\begin{equation}\label{eq:jrp-bilinear}
    \min_{X\in\Delta^1_{1,3}}\ \max_{\tY\in\Delta^K_{1,2},\ \tZ\in\Delta^{K_0}_2}\
    L(X,\tY,\tZ) := \langle H + \tY + \tZ,\ X\rangle ,
\end{equation}
which is \eqref{eq:saddle} with $\cX = \Delta^1_{1,3}$, $\cY = \Delta^K_{1,2}\times\Delta^{K_0}_2$, $c = \vc(H)$, $q = 0$ and $A : X\mapsto(-X,-X)$.
Normalized duals for \eqref{eq:jrp-lp} are recovered as $Y = \tY\oslash K$ and $Z = \tZ/K_0$.
\end{prop}

\begin{proof}
We prove Proposition~\ref{prop:jrp-reform} in two steps: price the covering rows, and absorb the simplices that pricing produces together with the assignment rows.

\emph{Step 1: the dual of each group is a simplex.}
Price $x_{ist} \le y_{is}$ with $\tY_{ist} \ge 0$ and $x_{ist} \le z_s$ with $\tZ_{ist} \ge 0$, so that the Lagrangian of \eqref{eq:jrp-lp} reads
\[
    L(X,y,z;\tY,\tZ) \;=\; \langle H,\ X\rangle
    \;+\; \sum_{i,s} y_{is}\Big(K_i - \sum_{t\ge s}\tY_{ist}\Big)
    \;+\; \sum_{s} z_s\Big(K_0 - \sum_{i}\sum_{t\ge s}\tZ_{ist}\Big)
    \;+\; \langle \tY + \tZ,\ X\rangle .
\]
Each $y_{is}$ and each $z_s$ occurs in no row but its own group's, and its nonnegativity is implied by that group together with $x \ge 0$, so it is free and the row its column contributes to the dual is an assignment row: $\sum_{t\ge s}\tY_{ist} = K_i$ for each item group $(i,s)$, and $\sum_i\sum_{t\ge s}\tZ_{ist} = K_0$ for each joint group $s$.

\emph{Step 2: absorb, on both sides.}
Definition~\ref{def:split} applied to those rows puts $\tY$ in $\Delta^K_{1,2}$ and $\tZ$ in $\Delta^{K_0}_2$, and applied to the assignment rows of \eqref{eq:jrp-lp} puts $X$ in $\Delta^1_{1,3}$.
Every constraint of \eqref{eq:jrp-lp} has been absorbed, on both sides: all structure now resides in three products of scaled simplices, and what is left is \eqref{eq:jrp-bilinear}.
\end{proof}

\paragraph{Step size.}
The groups of both dual families partition the primal index set, so $\tY$ and $\tZ$ are tensors of the same shape as $X$ and the coupling of \eqref{eq:jrp-bilinear} is the negated duplication $A : X\mapsto(-X,-X)$.
At $\Lambda = I$,
\begin{equation}\label{eq:jrp-Anorm}
    A^\top A = 2I,
    \qquad
    \cnorm{A} = \norm{A}_2 = \sqrt 2 ,
\end{equation}
so any $\eta < 1/\sqrt2$ is admissible at every primal weight.
Each $X_{ist}$ lies in exactly one item group $(i,s)$ and exactly one joint group $s$, so a diagonal rescaling would be the scalar $2I$, a change of step size.

\paragraph{Per-iteration cost.}
With $\Lambda = I$, $c = \vc(H)$, $q = 0$ and $A^\top(\tY,\tZ) = -\tY-\tZ$, the iteration \eqref{eq:pdhg} in these coordinates is
\begin{equation}\label{eq:jrp-pdhg}
\begin{aligned}
    X^{k+1} &:= \proj_{\Delta^1_{1,3}}\!\pran{X^k - \tfrac{\eta}{\omega}\pran{H + \tY^k + \tZ^k}},\\
    \tY^{k+1} &:= \proj_{\Delta^K_{1,2}}\!\pran{\tY^k + \eta\omega\pran{2X^{k+1}-X^k}},
    \qquad
    \tZ^{k+1} := \proj_{\Delta^{K_0}_2}\!\pran{\tZ^k + \eta\omega\pran{2X^{k+1}-X^k}},
\end{aligned}
\end{equation}
which is elementwise: forming $H+\tY+\tZ$ and the reflected term $2X^{k+1}-X^k$ is $O(NT^2)$ work.
The three projections decompose into independent simplex projections: $NT$ segments of length at most $T$ for $\Delta^1_{1,3}$ and $\Delta^K_{1,2}$, and $T$ segments of length up to $NT$ for $\Delta^{K_0}_2$.
Projecting each segment by median selection
\cite{brucker_on_1984,pardalos_algorithm_1990} is linear in its length, so a step is
$O(NT^2)$ work.
Sort-and-threshold \cite{held_validation_1974} costs a logarithmic factor more.
The implementation uses active-set methods
\cite{michelot_finite_1986,perez_efficient_2022}, which trade the worst-case guarantee
for data-dependent control flow that is faster in practice.
PDHG applied directly to \eqref{eq:jrp-lp} multiplies twice per iteration by a constraint matrix with roughly $\tfrac52NT^2$ nonzeros, with the irregular memory access of general sparse matrices.

\paragraph{Termination.}
Both iterates are feasible at every step, so the gap is computable in closed form:
\begin{equation}\label{eq:jrp-bounds}
    L_D(\tY,\tZ) = \sum_{i,t}\min_{s\le t}\,(H+\tY+\tZ)_{ist},
    \qquad
    L_P(X) = \langle H,X\rangle + \sum_{i,s}K_i\max_{t\ge s}X_{ist}
      + \sum_s K_0\max_{i,\,t\ge s}X_{ist},
\end{equation}
the first placing all mass on the cheapest source period, the second recovering the covering costs of the implied solution of \eqref{eq:jrp-lp}.
Both are exact objective values for \eqref{eq:jrp-lp}, so the relative gap the solver monitors is exact.

\subsubsection{Theoretical improvements}\label{sec:jrp-inst}

\begin{thm}[Convergence of PRIA]\label{thm:jrp-rate}
For the joint replenishment relaxation \eqref{eq:jrp-lp}, run the method of
Section~\ref{sec:alg} on \eqref{eq:jrp-bilinear} at a fixed $\omega = 1$ and
$0 < \eta < 1/\sqrt2$, on the normalized instance
$(\vc(H);K;K_0)/\norm{(\vc(H);K;K_0)}_2$.
Write $v^k$ for its iterates without restarts, $v^{j,k}$ for those of its $j$th
epoch under restarts, and $\norm{\cdot}_P$ for the canonical norm
\eqref{eq:canonical} at that primal weight.
\begin{enumerate}[label=(\alph*)]
    \item \textbf{Sublinear convergence without restarts.}
    For every $k \ge 0$,
    \[
        \norm{v^k - T_\omega(v^k)}_P \;\le\; \frac{2}{k+1}\sqrt{\frac{T(N+1)}{\eta}} .
    \]
    \item \textbf{Linear convergence with adaptive restarts.}
    Restarting by sufficient decay \eqref{eq:adaptive-restart}, every epoch has length at most
    $\lceil e/\alpha\rceil$ for the sharpness modulus $\alpha > 0$ of
    Definition~\ref{def:alpha}, and for every $j \ge 0$,
    \[
        \dist_P(v^{j+1,0},\cZ^\star) \;\le\;
        \frac{e^{-j}}{\alpha}\sqrt{\frac{T(N+1)}{\eta}} .
    \]
\end{enumerate}
\end{thm}

\begin{proof}
Every block is a product of scaled simplices, so $\Lambda = I$ and the reformulation has
\begin{equation}\label{eq:jrp-constants}
    \cnorm{A} \;=\; \sqrt2,
    \qquad
    \diam_\Lambda(\cZ)^2 \;\le\; 2NT + 2T\norm{K}_2^2 + 2TK_0^2 \;\le\; 2NT + 2T
\end{equation}
by \eqref{eq:jrp-Anorm} with $A^\top A = 2I$, and by
$\diam_\Lambda\pran{\prod_g\Delta^{r_g}}^2 \le 2\sum_g r_g^2$ over $NT$ unit simplices and $T$
blocks of each total $K_i$ and $K_0$.
Both blocks are polytopes, so $\alpha > 0$.
The initial point lies in $\cZ$, so \eqref{eq:lammax} at $\omega = 1$ gives
\[
    \dist_P(v^{0,0},\cZ^\star) \;<\; \sqrt{\frac{2}{\eta}}\sqrt{2T(N+1)}
    \;=\; 2\sqrt{\frac{T(N+1)}{\eta}} .
\]
Part (a) is Proposition~\ref{prop:luyang}(a) with this bound, part (b)
Proposition~\ref{prop:luyang}(b) with this $\alpha$ and $\tau^0 \ge 1$.
\end{proof}

\begin{table}[ht]
    \centering
    \small
    \begin{tabular}{@{}l l l l@{}}
        \toprule
        & Sublinear rate & Linear rate & Iteration cost \\
        \midrule
        PDLP formulation & $O\!\pran{\cnorm{A}R/\varepsilon}$
          & $O\!\pran{\cnorm{A}H_{\mathrm{KKT}}\sqrt{1+4\rho^2}\log(1/\varepsilon)}$ & $O(NT^2)$ \\
        Reformulation & $O\!\pran{\sqrt{NT}/\varepsilon}$
          & $O\!\pran{\alpha^{-1}\log(1/\varepsilon)}$ & $O(NT^2)$ \\
        \bottomrule
    \end{tabular}
\caption{Applications of $T_\omega$ to reach accuracy $\varepsilon$, over $N$ items and
$T$ periods, with arithmetic per iteration in the last column.
The columns are those of Table~\ref{tab:qubo-rates}.
Here $\alpha$ is the sharpness modulus of Definition~\ref{def:alpha}.}
    \label{tab:jrp-rates}
\end{table}

Table~\ref{tab:jrp-rates} presents the effect of reformulation on the constants governing convergence.
For the standard PDLP formulation, the sublinear bound involves the full constraint matrix through $\cnorm{A}$ and the initial distance $R$. Every entry of that matrix is $\pm1$, so $\cnorm{A}$ is fixed by $N$ and $T$. The distance $R$ depends on the costs through the saddle set. That set is known only once the instance is solved, so the bound stays loose.
For the reformulation, the coupling is the negated duplication $A:X\mapsto(-X,-X)$, so $\cnorm{A}=\sqrt2$ exactly, on every instance and at every primal weight. The three products of scaled simplices are bounded, which yields an explicit bound on the initial distance in terms of $N$ and $T$, providing an explicit and tight sublinear rate. Interestingly, both constants depend on the instance shape alone, so this sublinear bound is a function of $N$ and $T$, independent of the costs $K$, $K_0$ and $H$.
Both rows obtain their linear rate from the same sharpness. For the original PDLP formulation it enters through $H_{\mathrm{KKT}}$ together with the radius $\rho$, which are data dependent and hard to characterize. For the reformulation it is the sharpness modulus $\alpha$ of Definition~\ref{def:alpha}, which is positive because both blocks are polytopes, so Theorem~\ref{thm:jrp-rate}(b) is a linear rate with a finite constant.
The problem-dependent terms in the sublinear rate are therefore $N$ and $T$ alone, and the linear rate measures its sharpness on a saddle problem over three products of simplices in place of the KKT system of \eqref{eq:jrp-lp}.
The different residuals and norms specified in the caption must be taken into account when comparing these bounds.

Reformulation also changes the operations within each iteration.
Although both formulations have $O(NT^2)$ arithmetic cost, the reformulation replaces the two products against the $\approx\tfrac52NT^2$ nonzeros of \eqref{eq:jrp-lp} with the negated duplication $A:X\mapsto(-X,-X)$, which is elementwise. This permits a GPU implementation that stores no constraint matrix at all.
The cost moves to the projections. Here $\cX$ and $\cY$ are products of scaled simplices, so each iteration projects onto $\Delta^1_{1,3}$, $\Delta^K_{1,2}$ and $\Delta^{K_0}_2$, where the standard form clips coordinates to a box.
The simplex projections dominate, and they make PRIA's iterations the slower of the two (Table~\ref{tab:jrp-cost}).
The two channels of Section~\ref{sec:whynum} therefore act in opposite directions on the JRP, and the running-time advantage comes through the iteration count.

\subsubsection{Computational results}\label{sec:jrp-results}

\paragraph{Implementation.}
The product $AX$ is computed by negating $X$ into the two dual tensors, and the transpose product $A^\top(\tY,\tZ)$ by negating their sum, both elementwise over tensors of one shape.
The primal projection onto $\Delta^1_{1,3}$ and the dual projections onto $\Delta^K_{1,2}$ and $\Delta^{K_0}_2$ are segmented simplex projections, sort-free by default, using Michelot-type active-set thresholding \cite{michelot_finite_1986} with warm-started thresholds and a compacted active set for the large per-period segments of $\tZ$.
The remaining iterate updates are elementwise tensor operations.
The upper-triangular tensors are packed densely: with $T$ padded to an even integer, pairing column $t$ with column $T-1-t$ makes every packed row exactly full, and all layout conversions are closed-form index maps implemented as gathers.

\paragraph{Instances.}
Given $N$ and $T$, a random instance is drawn with integer costs.
Order costs are $N+1$ i.i.d.\ values $u\sim\mathrm{Uniform}(1,T)$, rounded to the nearest
integer with a floor of $1$.
The first $N$ are the $K_i$ and the last is $K_0$.
Sampling $K_0$ on the order of the horizon places instances near the empirically hardest
regime for the relaxation.
Holding costs: a slope $a_i\sim\mathrm{Uniform}(0,1)$ per item and an offset $\varepsilon_{ist}\sim\mathrm{Uniform}(0,a_i)$ per pair $s<t$, with $H_{ist} = \mathrm{round}(a_i(t-s)+\varepsilon_{ist})$ and $H_{ist} = 0$ for $s \ge t$, so holding costs are nonnegative and increasing in the holding interval.
\paragraph{Protocol.}
We test the eleven shapes of Table~\ref{tab:instance-sizes}, with more repetitions at the cheaper sizes, at a relative optimality tolerance $10^{-6}$, relative feasibility tolerance $10^{-9}$ and a five-minute limit per instance.
The baselines are cuPDLPx \cite{lu_cupdlpx_2025} and Gurobi \cite{gurobi_optimization_llc_gurobi_2026}, each solving \eqref{eq:jrp-lp}.
PRIA and cuPDLPx each ran on one NVIDIA H100 80\,GB.
Gurobi ran on one $64$-core CPU node.

\begin{table}[ht]
    \centering
    \begin{tabular}{llrrr}
        \toprule
        Family & $N \times T$ & $x$ variables & Nonzeros in \eqref{eq:jrp-lp} & Instances \\
        \midrule
        \multirow{4}{*}{\texttt{medium}}
          & $10\times25$   & $3{,}250$      & $\approx 0.02$M & $100$ \\
          & $20\times50$   & $25{,}500$     & $\approx 0.13$M & $75$ \\
          & $30\times75$   & $85{,}500$     & $\approx 0.4$M  & $60$ \\
          & $50\times100$  & $252{,}500$    & $\approx 1.3$M  & $45$ \\
        \midrule
        \multirow{4}{*}{\texttt{large}}
          & $75\times150$  & $849{,}375$    & $\approx 4$M    & $30$ \\
          & $100\times200$ & $2{,}010{,}000$& $\approx 10$M   & $20$ \\
          & $150\times250$ & $4{,}706{,}250$& $\approx 23$M   & $12$ \\
          & $200\times300$ & $9{,}030{,}000$& $\approx 45$M   & $8$ \\
        \midrule
        \multirow{3}{*}{\texttt{huge}}
          & $250\times375$ & $17{,}625{,}000$& $\approx 88$M  & $6$ \\
          & $300\times400$ & $24{,}060{,}000$& $\approx 120$M & $4$ \\
          & $350\times450$ & $35{,}516{,}250$& $\approx 177$M & $3$ \\
        \bottomrule
    \end{tabular}
\caption{The benchmark ladder: eleven shapes over four decades of $x$ variables, $363$ instances in all.
Seed counts taper with shape, since the cost of a shape grows faster than the information another seed at it buys.
Nonzero counts are $\approx\tfrac52NT^2$ for \eqref{eq:jrp-lp}.}
    \label{tab:instance-sizes}
\end{table}

\paragraph{Results.}
Table~\ref{tab:jrp-solve} is the running time comparison.
PRIA gains on cuPDLPx and Gurobi in solve time at every instance size, by factors of $6.1$--$28.9$ over cuPDLPx.
PRIA solved all $363$ instances, where cuPDLPx solved $347$ and Gurobi $308$, so the entries above $10^5$ variables are floors and the margins there are lower bounds.

\begin{table}[ht]
    \centering
    \small
    \begin{tabular}{@{}l r r rr rr r@{}}
        \toprule
        & & & \multicolumn{2}{c}{cuPDLPx} & \multicolumn{2}{c}{Gurobi} & \\
        \cmidrule(lr){4-5}\cmidrule(lr){6-7}
        $x$ variables & Instances & PRIA & s & solved & s & solved & cuPDLPx/PRIA \\
        \midrule
        $10^3$--$10^4$   & 100 & \textbf{0.031} & 0.187   & 100 & 0.087   & 100 & $6.1$ \\
        $10^4$--$10^5$   & 135 & \textbf{0.136} & 1.32    & 135 & 2.14    & 135 & $9.7$ \\
        $10^5$--$10^6$   & 75  & \textbf{0.61}  & $17.5$ & 74  & $72.9$ & 73  & $28.9$ \\
        $10^6$--$10^7$   & 40  & \textbf{2.88}  & $80.4$ & 30  & ---     & 0   & $27.9$ \\
        $10^7$--$10^8$   & 13  & \textbf{11.3}  & $234$  & 8   & ---     & 0   & $20.8$ \\
        \bottomrule
    \end{tabular}
\caption{JRP relaxation: solve time in seconds, shifted (by $10^{-2}$\,s) geometric mean, by number of $x$ variables.
A paired column gives how many of the bucket finished within the five-minute limit; an unfinished bucket's entry is a floor. PRIA solved every instance.}
    \label{tab:jrp-solve}
\end{table}

Tables~\ref{tab:jrp-cost} and~\ref{tab:jrp-iters} show that the running-time advantage of PRIA comes from the iteration count alone.
Across the five size groups, PRIA takes $1.3$--$2.9$ times as long per iteration as cuPDLPx and converges in $13.6$--$52.3$ times fewer iterations.
The higher iteration cost is consistent with projecting onto products of scaled simplices where the standard form clips coordinates, and in particular with the $\Delta^{K_0}_2$ block, whose $T$ segments run to length $NT$.
The two mechanisms discussed in Section~\ref{sec:whynum} therefore act against each other here, and the iteration count decides the comparison.

Table~\ref{tab:jrp-memory} shows a corresponding reduction in peak device memory.
For the two largest size groups, cuPDLPx uses $3.5$--$3.8$ times as much device memory as PRIA.
In particular, for instances with $10^7$--$10^8$ variables, PRIA uses $5.29$\,GiB compared with $20.1$\,GiB for cuPDLPx.
This is consistent with representing the coupling as a duplication and storing no constraint matrix, where \eqref{eq:jrp-lp} carries $\approx\tfrac52NT^2$ nonzeros.

\begin{table}[ht]
    \centering
    \small
    \begin{minipage}[t]{0.48\linewidth}
    \centering
    \begin{tabular}{@{}l rr r@{}}
        \toprule
        $x$ variables & PRIA & cuPDLPx & cuPDLPx/PRIA \\
        \midrule
        $10^3$--$10^4$ & 123 & \textbf{41.3} & $0.34$ \\
        $10^4$--$10^5$ & 59.4 & \textbf{32.4} & $0.55$ \\
        $10^5$--$10^6$ & 131 & \textbf{71.9} & $0.55$ \\
        $10^6$--$10^7$ & 825 & \textbf{583} & $0.71$ \\
        $10^7$--$10^8$ & 4960 & \textbf{3890} & $0.78$ \\
        \bottomrule
    \end{tabular}
\caption{Median microseconds per iteration.
PRIA's iterations are slower because of the simplex projections, and in particular the $\Delta^{K_0}_2$ block, whose $T$ segments run to length $NT$.}
    \label{tab:jrp-cost}
    \end{minipage}\hfill
    \begin{minipage}[t]{0.48\linewidth}
    \centering
    \begin{tabular}{@{}l rr r@{}}
        \toprule
        $x$ variables & PRIA & cuPDLPx & cuPDLPx/PRIA \\
        \midrule
        $10^3$--$10^4$ & \textbf{245} & 3{,}338 & $13.6$ \\
        $10^4$--$10^5$ & \textbf{2{,}318} & 37{,}914 & $16.4$ \\
        $10^5$--$10^6$ & \textbf{3{,}611} & 189{,}002 & $52.3$ \\
        $10^6$--$10^7$ & \textbf{3{,}317} & 130{,}281 & $39.3$ \\
        $10^7$--$10^8$ & \textbf{2{,}363} & 61{,}329 & $26.0$ \\
        \bottomrule
    \end{tabular}
\caption{Iterations to converge, shifted geometric mean with shift $10$. Instances where cuPDLPx failed to converge were not included.}
    \label{tab:jrp-iters}
    \end{minipage}
\end{table}

\begin{table}[ht]
    \centering
    \small
    \begin{tabular}{@{}l r rr r@{}}
        \toprule
        $x$ variables & Instances & PRIA & cuPDLPx & cuPDLPx/PRIA \\
        \midrule
        $10^6$--$10^7$ & 40 & \textbf{1.60\,GiB} & 5.60\,GiB & $3.5$ \\
        $10^7$--$10^8$ & 13 & \textbf{5.29\,GiB} & 20.1\,GiB & $3.8$ \\
        \bottomrule
    \end{tabular}
\caption{Largest device memory held during a solve, sampled with NVML.}
    \label{tab:jrp-memory}
\end{table}

\section{Conclusions}\label{sec:conclusions}

Reformulation, long considered of limited value for interior-point methods, is a substantial lever for PDHG-type LP solvers, and it acts through two independent channels.
On both of our instances, absorbing every constraint into products of projection-friendly polyhedra removed the sparse constraint matrix from the iteration and changed the geometry of the saddle problem.
Both reformulated solvers run faster than cuPDLPx on the original formulations at every rung of both ladders (Tables~\ref{tab:qubo-solve} and~\ref{tab:jrp-solve}).

We study what to expect from a new instantiation and what to check.
The coupling norm is a property of the split, made data-free by the group structure it
absorbs.
The diameter is computed from the absorbed sets once they are known to be bounded.
Sharpness is finite for polyhedral blocks, and bounded when the reformulated polytope
carries a lattice of vertex values.

Several directions follow.
The recipe applies wherever auxiliary variables select maxima over groups, beyond inventory and binary quadratic models: assignment blocks, $\ell_1$ structures, and max-selection variables appear across large-scale LP applications, and each is an absorption candidate by Table~\ref{tab:constructs}.
A presolve pass that \emph{detects} embedded projection-friendly blocks and hands a general-purpose PDHG solver the reformulated problem would bring these gains without modelling effort.
On the theory side, the restriction to polyhedral blocks and the restriction to bounded ones are likely both removable, at the cost of a normalized progress measure and a H\"olderian error bound.

\section*{Acknowledgments}
We are grateful for the very helpful discussion with Jinwen Yang on the JRP reformulation scheme that initiated this line of research. We thank Edward Rothberg for suggesting, during a visit to MIT, that QUBO would be a good candidate for reformulation. We thank Zedong Peng for the discussions about cuPDLPx.
Haihao Lu thanks Oliver Hinder, Warren Schudy, Miles Lubin and David Applegate for their early discussions about reformulation for PDLP.

Jacob Dentes is supported by the National Science Foundation Graduate Research Fellowship Program under Grant No.\ 2141064.
Any opinions, findings, and conclusions or recommendations expressed in this material are those of the author(s) and do not necessarily reflect the views of the National Science Foundation.
Haihao Lu is partially supported by AFOSR Grant No.\ FA9550-24-1-0051, ONR Grant No.\ N000142412735 and a Sloan Research Fellowship.

\bibliography{ref-papers}

\appendix

\section{Sharpness of the QUBO reformulation}\label{app:sharp}

This appendix proves \eqref{eq:qubo-alpha}, the sharpness constant of
Theorem~\ref{thm:qubo-rate}.
The rates of that theorem are Proposition~\ref{prop:luyang} applied to
\eqref{eq:qubo-constants} and \eqref{eq:qubo-alpha}.

\begin{prop}[Sharpness of the QUBO reformulation]\label{prop:qubo-sharp}
Let $W$ satisfy Assumption~\ref{ass:qubo-integral}, let every node carry at least one edge, and let $n \ge 5$.
Run the method at a fixed $\omega = 1$ and $0 < \eta < 1$ on the normalized instance
$(d;w)/\norm{(d;w)}_2$, and write $\norm{\cdot}_P$ for the canonical norm \eqref{eq:canonical} at that
primal weight.
Then
\begin{equation}\label{eq:sharp-final}
    \dist_P(v,\cZ^\star) \;\le\;
    \pran{1 + \frac{16}{\eta}\norm{W}_F\,m\sqrt n}\norm{v - T_\omega(v)}_P
    \qquad\text{for every } v \in \RR^{n+m},
\end{equation}
where $\norm{W}_F$ is the Frobenius norm of $W$ itself.
\end{prop}

Every node must carry an edge, or $\Lambda_{\cX}$ and $P$ are singular and $T_\omega$ is
undefined.
This costs nothing: a node with no edge appears in no quadratic term, so the sign of its
linear coefficient fixes its value.
Section~\ref{sub:margins} bounds $\gamma$ by the vertex margins of two auxiliary LPs.
Section~\ref{sub:halfint} shows that both LPs have a half-integral lattice.
Section~\ref{sub:residual} converts that bound into \eqref{eq:sharp-final}.

\subsection{Vertex margins}\label{sub:margins}

Write $\gamma$ for the smallest constant with $\dist_\Lambda(v,\cZ^\star) \le \gamma\,\gp(v)$ on
$\cZ$.
Lemma~\ref{lem:sharp} turns a bound on $\gamma$ into one on $1/\alpha$, and the rest of
this appendix bounds $\gamma$.
The gap between vertex values controls how fast a linear objective grows away from its
optimal face.

\begin{lem}[Vertex margins give sharpness]\label{lem:vertexmargin}
Let $\cQ \subseteq \RR^N$ be a polytope, $g \in \RR^N$, $\mathrm{OPT} := \min_{u\in\cQ}\langle g,u\rangle$,
$\cF$ the optimal face, and
\begin{equation}\label{eq:margin}
    \delta \;:=\; \min\left\{\langle g,\bar u\rangle - \mathrm{OPT} \;:\;
    \bar u \text{ a vertex of } \cQ \text{ with } \langle g,\bar u\rangle > \mathrm{OPT}\right\}
\end{equation}
(with $\delta := +\infty$ if every vertex is optimal).
Then, for $\diam_2(\cQ) > 0$,
\[
    \langle g, u\rangle - \mathrm{OPT} \;\ge\; \frac{\delta}{\diam_2(\cQ)}\,\dist_2(u,\cF)
    \qquad\text{for every } u \in \cQ .
\]
\end{lem}

\begin{proof}
Write $u = \sum_{\bar u}\mu_{\bar u}\bar u$ over the vertices of $\cQ$.
Put $p := \sum_{\bar u \in \cF}\mu_{\bar u}$.
If $p = 1$ then $u \in \cF$ and both sides vanish, so assume $p < 1$.
Set $b := (1-p)^{-1}\sum_{\bar u\notin\cF}\mu_{\bar u}\bar u \in \cQ$.
Set $a := p^{-1}\sum_{\bar u\in\cF}\mu_{\bar u}\bar u \in \cF$ if $p > 0$, and let $a$ be
any point of $\cF$ if $p = 0$.
Then $u = pa+(1-p)b$, and
\[
\begin{aligned}
    \dist_2(u,\cF) &\;\le\; \norm{u-a} \;=\; (1-p)\norm{b-a} \;\le\; (1-p)\diam_2(\cQ), \\
    \langle g,u\rangle - \mathrm{OPT}
      &\;=\; \sum_{\bar u}\mu_{\bar u}\pran{\langle g,\bar u\rangle - \mathrm{OPT}}
      \;\ge\; (1-p)\delta .
\end{aligned}
\]
The first display uses $a \in \cF$ and $a,b \in \cQ$.
The second uses linearity of $\langle g,\cdot\rangle$.
The vertices in $\cF$ are minimizers and contribute zero; those outside contribute at
least $\delta$ each.
Dividing the second display by the first gives the claim.
\end{proof}

$L_P$ and $L_D$ are the value functions of two linear programs.
The \emph{primal value LP} is \eqref{eq:qubo-lp} itself.
Its feasible set $\cQ_P$ is cut out by the covering and linking rows, $x\in[0,1]^n$ and
$z\in[0,1]^E$, and minimizing $z$ out of it gives $L_P$.
The \emph{dual value LP} is the hypograph form of $L_D$ over the absorbed box $\cY$.
Lemma~\ref{lem:qubo-dualmargin} writes it out.
Write $g_P,g_D$ for their objective vectors, $D_P,D_D$ for the Euclidean diameters of
$\cQ_P,\cQ_D$, and $\delta_P,\delta_D$ for their vertex margins \eqref{eq:margin}.
Take the dual margin on $-g_D$, because that LP is a maximization.

\begin{lem}[Value polytopes give sharpness]\label{lem:gammavertex}
If $\cQ_P$ and $\cQ_D$ are bounded then
\[
    \gamma \;\le\; \max\left\{\norm{\Lambda_{\cX}}_2^{1/2}\frac{D_P}{\delta_P},\
                              \norm{\Lambda_{\cY}}_2^{1/2}\frac{D_D}{\delta_D}\right\},
\]
with all four quantities measured on the same instance.
\end{lem}

\begin{proof}
Fix $x\in\cX$ and lift it to $u\in\cQ_P$ by restoring the eliminated variables at their
optimal values, so that $\langle g_P,u\rangle = L_P(x)$.
Write $\pi$ for the projection onto the $x$ coordinates.
Then $\pi$ carries $\cF_P$ into $\cX^\star$: if $(x',z')\in\cF_P$ then
$L_P(x') \le \langle g_P,(x',z')\rangle = \mathrm{OPT} \le L_P(x')$.
For every $f\in\cF_P$, $1$-Lipschitzness gives
$\dist_2(x,\cX^\star) \le \norm{\pi u - \pi f} \le \norm{u-f}$.
Taking the infimum over $f$ gives $\dist_2(x,\cX^\star) \le \dist_2(u,\cF_P)$.
Lemma~\ref{lem:vertexmargin} then gives
$L_P(x) - \mathrm{OPT} \ge (\delta_P/D_P)\dist_2(x,\cX^\star)$.
The dual value LP is a maximization, so the same argument on $-g_D$ gives
$\mathrm{OPT} - L_D(\lambda) \ge (\delta_D/D_D)\dist_2(\lambda,\cY^\star)$.

Write $s := (\delta_P/D_P)\dist_2(x,\cX^\star)$ and
$t := (\delta_D/D_D)\dist_2(\lambda,\cY^\star)$.
Write $C_P,C_D$ for the squares of the two terms in the stated maximum. Then
\[
\begin{aligned}
    \gp(v) &\;=\; \pran{L_P(x)-\mathrm{OPT}} + \pran{\mathrm{OPT}-L_D(\lambda)} \;\ge\; s+t, \\
    \dist_\Lambda(v,\cZ^\star)^2 &\;\le\; C_Ps^2 + C_Dt^2 \;\le\; \max\{C_P,C_D\}(s+t)^2 .
\end{aligned}
\]
The first line adds the two bounds.
The second uses that $\cZ^\star = \cX^\star\times\cY^\star$ is a product and $\Lambda$ is
block diagonal, so that
\[
\begin{aligned}
    \dist_\Lambda(v,\cZ^\star)^2
      &\;=\; \dist_{\Lambda_{\cX}}(x,\cX^\star)^2 + \dist_{\Lambda_{\cY}}(\lambda,\cY^\star)^2 , \\
    \dist_{\Lambda_{\cX}}(x,\cX^\star)^2
      &\;\le\; \norm{\Lambda_{\cX}}_2\,\dist_2(x,\cX^\star)^2 \;=\; C_Ps^2 , \\
    \dist_{\Lambda_{\cY}}(\lambda,\cY^\star)^2
      &\;\le\; \norm{\Lambda_{\cY}}_2\,\dist_2(\lambda,\cY^\star)^2 \;=\; C_Dt^2 ,
\end{aligned}
\]
and then adds the nonnegative cross term $2\max\{C_P,C_D\}st$.
Taking square roots gives
$\dist_\Lambda(v,\cZ^\star) \le \sqrt{\max\{C_P,C_D\}}\,\gp(v)$, whose constant is the
stated maximum.
\end{proof}

One hypothesis supplies both margins.

\begin{cor}[A lattice of vertices makes $\gamma$ explicit]\label{cor:lattice}
Suppose the objective data of both value LPs is integral and every vertex of $\cQ_P$ and
of $\cQ_D$ has all coordinates in $\tfrac1q\ZZ$ for an integer $q \ge 1$.
Then $\delta_P,\delta_D \ge 1/q$.
\end{cor}

\begin{proof}
An integral linear functional evaluated at a point of $\tfrac1q\ZZ^N$ lies in
$\tfrac1q\ZZ$, and $\mathrm{OPT}$ is attained at a vertex of a bounded polytope, so it is
itself such a value; distinct values of $\tfrac1q\ZZ$ differ by at least $1/q$.
\end{proof}

\subsection{Both value polytopes are half-integral}\label{sub:halfint}

\begin{lem}[Every vertex is half-integral]\label{lem:qubo-halfint}
Every vertex of $\cQ_P$ satisfies $(x,z) \in \{0,\tfrac12,1\}^{n+m}$.
\end{lem}

\begin{proof}
Let $(x,z)$ be a vertex of $\cQ_P$.
Let $\cB$ be a set of $n+m$ independent tight rows.
No edge is a self-loop, because $E\subseteq\{i<j\}$.
A self-loop would put a coefficient $2$ into a residual equation.

The variable $z_e$ occurs only in the rows of its own edge.
If $\cB$ held none of them, the column of $z_e$ would vanish.
So $\cB$ holds $t_e \ge 1$ rows of each edge.
Let $p$ be the number of rows of $\cB$ that bound $x$. Then
\[
    \sum_e t_e + p \;=\; n+m .
\]
Choose one row per edge.
Its $z_e$ coefficient is $\pm1$.
It therefore fixes $z_e \in \{0,\ 1,\ x_i,\ x_j,\ x_i+x_j-1\}$.
Order the columns $(z,x)$ and the rows (chosen, rest).
The chosen block is then diagonal with entries $\pm1$.
Writing $\cB$ also for the matrix of those rows, $\det\cB = \pm\det S$ for the Schur
complement $S$, the system left on $x$.
$S$ has $\sum_e(t_e-1)+p = n$ rows and rank $n$, for any admissible choice of rows.

Each remaining row of edge $e$ combines with the chosen one and cancels $z_e$.
A covering row has three variables.
Every such combination has at most two entries of $x$:
\[
\begin{aligned}
    \Ep:\ &\{\text{cover},\ z_e\ge0\} \Rightarrow x_i+x_j = 1,
          &&\{\text{cover},\ z_e\le1\} \Rightarrow x_i+x_j = 2, \\
    \Em:\ &\{z_e\le x_i,\ z_e\le x_j\} \Rightarrow x_i = x_j,
          &&\{z_e\le x_i,\ z_e\ge0\} \Rightarrow x_i = 0, \\
          &\{z_e\le x_i,\ z_e\le1\} \Rightarrow x_i = 1 .
\end{aligned}
\]
The pair $\{z_e\ge0,\ z_e\le1\}$ is infeasible. The rows at $x_j$ are symmetric.
With the $p$ bound rows, every row of $S$ reads
\begin{equation}\label{eq:qubo-switch}
    x_j = \epsilon x_i + c \quad(\epsilon\in\{\pm1\},\ c\in\ZZ),
    \qquad\text{or}\qquad
    x_i = c \quad (c \in \{0,1\}).
\end{equation}
Each edge gives at most one equation of the first kind.
Distinct edges join distinct pairs.
Those equations therefore form a simple graph on $[n]$.

Fix a component $C$ of that graph and a root $r$, and set $x_r = t$.
Each equation flips a sign and shifts by an integer.
A spanning tree therefore gives $x_v = \epsilon_vt+c_v$ with $\epsilon_v\in\{\pm1\}$ and
$c_v\in\ZZ$.
Every row of $S$ has support in one component, so $S$ is block diagonal.
No block has fewer than $|C|$ rows, or it would be rank deficient.
The rows number $n = \sum_C|C|$.
Each block therefore has exactly $|C|$ rows, one more than its spanning tree.
The extra row is a pin or a non-tree equation.
A pin $x_v = c$ gives $t = \epsilon_v(c-c_v)\in\ZZ$.
A non-tree equation must disagree in sign with the tree path.
The tree rows have kernel spanned by $(\epsilon_v)_{v\in C}$.
A sign-consistent row vanishes on that kernel and leaves the block singular.
Disagreement gives $\epsilon t + c = -\epsilon t + c'$, hence $2t\in\ZZ$.
In both cases $t \in \tfrac12\ZZ$, so $x \in \tfrac12\ZZ^n$.
With $0\le x\le1$ this gives $x \in \{0,\tfrac12,1\}^n$.
Each $z_e$ is one of $0,1,x_i,x_j,x_i+x_j-1$.
With $z\in[0,1]^E$ this gives $z_e \in \tfrac12\ZZ\cap[0,1]$.
\end{proof}

The classical statement \cite{hammer_roof_1984} is about optima.
Lemma~\ref{lem:qubo-halfint} extends it to all vertices, which is what
Corollary~\ref{cor:lattice} needs.
The dual side needs a lemma on systems in which every variable touches at most two
equations.
Such a system can have an exponentially large determinant.
The lemma bounds the \emph{solution} instead, because the blocks are solved
independently.

\begin{lem}[Two-nonzero systems]\label{lem:twononzero}
Let $A$ be a nonsingular square matrix with entries in $\{0,\pm1\}$ having at most two
nonzero entries in every column.
Then $A^{-1}r \in \tfrac12\ZZ^N$ for every integral $r$.
\end{lem}

\begin{proof}
Form the multigraph $G$ on the rows of $A$ together with one extra vertex $\infty$.
Its edges are the columns: a column with two nonzeros joins its two rows, and a column
with one nonzero joins its row to $\infty$.
A column with no nonzero would make $A$ singular.
Let $C_1,\dots,C_k$ be the components of $G-\infty$.
Every column has its row-nonzeros inside one component.
After permuting, $A = \diag(A_1,\dots,A_k)$.
Each $A_t$ is square, since a block with more columns than rows would have a kernel, and
each is nonsingular.
Cramer's rule inside a block gives $(A^{-1}r)|_{C_t} = A_t^{-1}r_{C_t}$.
Its numerators are determinants of integral matrices and its denominators divide
$\det A_t$.
It is enough to prove $|\det A_t| \in \{1,2\}$.
The signs within a column are arbitrary, so $A_t$ shares the support of an oriented
incidence matrix while its signs stay free.
This matters below: an oriented incidence matrix carries one $+1$ and one $-1$ per
column, so its columns sum to zero, $\mathbf 1$ is a left null vector of any block with
no row deleted, and the second case would be vacuous.

Fix a block with $n_t$ rows and columns.
Let $H$ be the subgraph of $G$ induced on $C_t\cup\{\infty\}$, with $n_t+1$ vertices and
$n_t$ edges.
Suppose $H$ is connected.
It has one fewer edge than vertices, so it is a tree, and $A_t$ has the support of its
incidence matrix with the row of $\infty$ deleted.
A tree with at least two vertices has at least two leaves, so it has one other than
$\infty$.
That leaf's row has a single nonzero.
Expanding along it deletes the leaf and its edge and leaves a smaller tree.
Induction gives $|\det A_t| = 1$ for any signs.

Now suppose $H$ is disconnected.
$C_t$ is connected, so the components are $C_t$ and the isolated $\infty$.
No column of the block has a single nonzero, and $C_t$ has $n_t$ vertices and $n_t$
edges.
While some row has exactly one nonzero, expand along it.
This deletes a degree-one vertex with its edge.
It preserves connectivity and the equality $\#V = \#E$.
The process ends with a connected graph in which every vertex has degree at least two and
$\#V = \#E$.
Summing degrees gives $2\#E \ge 2\#V = 2\#E$, so every vertex has degree exactly two.
The graph is a single cycle.
In the determinant expansion a permutation matches each vertex to an incident edge and
uses each edge once.
Around a cycle the only two such matchings are the two orientations, so
\[
    \det A_t \;=\; \pm\prod_l s_l \;\pm\; \prod_l t_l ,
    \qquad s_l,t_l\in\{\pm1\} \text{ the two entries of column } l .
\]
Hence $|\det A_t| \in \{0,2\}$, and nonsingularity leaves $2$.
\end{proof}

\begin{lem}[Dual margin]\label{lem:qubo-dualmargin}
Under Assumption~\ref{ass:qubo-integral}, every vertex of $\cQ_D$ has all coordinates in $\tfrac12\ZZ$.
\end{lem}

\begin{proof}
The absorbed dual on $W$'s own data is one scalar $\lambda_e\in[0,w_e]$ per edge.
Write $a_e := \lambda_e$ on $\Em$ and $b_e := \lambda_e$ on $\Ep$, so $(a,b)$ has $m$
coordinates.
With
\[
    g_i(a,b) = d_i - \sum_{e\in\Em,\ i=\min e}a_e - \sum_{e\in\Em,\ i=\max e}(w_e-a_e)
             + \sum_{e\in\Ep,\ e\ni i}b_e ,
\]
the dual value LP is $\max\{\sum_i\varphi_i - \sum_{e\in\Ep}b_e\}$ over
\[
    \cQ_D = \left\{(a,b,\varphi) \;:\; 0\le a\le w,\ 0\le b\le w,\ -B\le\varphi\le0,\
    \varphi_i \le g_i(a,b)\right\},
    \qquad B_i := |d_i| + \sum_{e\ni i}w_e .
\]
The bound $\varphi\ge-B$ truncates $\cQ_D$ to a polytope.
At any optimum $\varphi_i = \min\{0,g_i(a,b)\}$, since raising $\varphi_i$ strictly improves
the objective.
On the box each $\Em$ term of $g_i$ lies in $[-w_e,0]$ and each $\Ep$ term in $[0,w_e]$,
so $|g_i(a,b)| \le B_i$.
The bound $\varphi_i \ge -B_i$ therefore holds at every optimum of the untruncated
problem.
The truncation leaves the optimal face intact and makes $\cQ_D$ bounded, as
Lemma~\ref{lem:gammavertex} requires.
Assumption~\ref{ass:qubo-integral} makes $B_i$ integral, which Lemma~\ref{lem:twononzero} needs for the
survivor $g_i(a,b) = -B_i$ below.

Let $(a,b,\varphi)$ be a vertex of $\cQ_D$.
Let $\cB$ be a set of $m+n$ independent tight rows.
The variable $\varphi_i$ occurs only in $\varphi_i\le0$, $\varphi_i\le g_i$ and
$\varphi_i\ge-B_i$, and no row involves two different $\varphi_i$.
Let $n_i$ count the rows of $\cB$ involving $\varphi_i$, and $p$ the tight bounds on
$(a,b)$. Then
\[
    \sum_i n_i + p \;=\; m+n .
\]
Here $n_i \ge 1$, else the column of $\varphi_i$ vanishes.
Also $n_i \le 2$: all three tight forces $0 = \varphi_i = -B_i$, and at $B_i = 0$ the rows
$\varphi_i\le0$ and $-\varphi_i\le0$ are negatives, so no independent $\cB$ holds both.
Spend one row per $i$ on $\varphi_i$, whose coefficient there is $\pm1$.
As in Lemma~\ref{lem:qubo-halfint}, the remaining $\sum_i(n_i-1)+p = m$ rows are
equations in $(a,b)$ of rank $m$.
For $n_i = 2$ the survivor combines the two tight rows,
\[
    \{\varphi_i\le0,\ \varphi_i\le g_i\} \;\Rightarrow\; g_i(a,b) = 0,
    \qquad
    \{\varphi_i\le g_i,\ \varphi_i\ge-B_i\} \;\Rightarrow\; g_i(a,b) = -B_i .
\]
The third pair is the dependent one just excluded.
Both right-hand sides are integral.
Writing a $g$-row with $(a,b)$ on the left moves
$-d_i + \sum_{e\in\Em,\ i=\max e}w_e$ to the right, which Assumption~\ref{ass:qubo-integral} makes integral too.
The $p$ bound rows pin variables to $0$ or $w_e$, both integers.
They pin $p$ distinct variables, because $w_e = 2|W_{ij}| > 0$ keeps $a_e\ge0$ and
$a_e\le w_e$ from both being tight.

Substituting the pins leaves $m-p$ equations in the $m-p$ unpinned variables.
Order the rows (pins, $g$-rows) and the columns (pinned, unpinned).
The system is then $\left(\begin{smallmatrix}\pm I_p&0\\ *&S\end{smallmatrix}\right)$, so
$\det S = \pm\det(\text{residual}) \ne 0$.
The reduced system is square and nonsingular, as Lemma~\ref{lem:twononzero} requires.
Its entries are $\{0,\pm1\}$ with at most two per column.
The reason is that $a_e$ and $b_e$ occur only in the two $g_i$ of the endpoints of $e$:
$a_e$ carries $-1$ at $\min e$ and $+1$ at $\max e$, and $b_e$ carries $+1$ at both.
Lemma~\ref{lem:twononzero} gives $(a,b) \in \tfrac12\ZZ^m$.
Hence $g_i(a,b)\in\tfrac12\ZZ$ and
$\varphi_i \in \{0, g_i(a,b), -B_i\} \subseteq \tfrac12\ZZ$.
\end{proof}

Both propositions are tight.
On the primal side a gap of exactly $\tfrac12$ occurs at $n = 3$.
On the dual side take
$W = \left(\begin{smallmatrix}-2&3&3\\3&-2&3\\3&3&-3\end{smallmatrix}\right)$, whose
residual system is the unsigned incidence matrix of the triangle, of determinant $-2$.
Its dual optimum $-\tfrac72$ is attained at one vertex, $b = (\tfrac12,\tfrac32,\tfrac32)$
and $\varphi = 0$, and the next vertex value is $-4$.

\subsection{From the gap to the fixed-point residual}\label{sub:residual}

The two lemmas of this subsection hold for any complete split; Section~\ref{sec:jrp-inst}
applies the second of them to the JRP.
Distances are measured in $\norm{\cdot}_\Lambda$, so residuals are measured in the dual
norm $\norm{\xi}_{\Lambda^{-1}} := \sqrt{\xi^\top\Lambda^{-1}\xi}$.
Write $\lambda_{\max}$ for the largest eigenvalue of $P$ relative to $\Lambda$, so that
$P \preceq \lambda_{\max}\Lambda$ and, by \eqref{eq:lammax} at $\omega = 1$,
$\lambda_{\max} < 2/\eta$.

\begin{lem}[Gap bounded by the residual]\label{lem:gapres}
For every $v \in \cZ$, $\gp(v) \le \diam_\Lambda(\cZ)\cdot\dist_{\Lambda^{-1}}(0,F(v))$, and for
every $v$ whatsoever, with $v^+ = T_\omega(v)$,
\[
    \dist_{\Lambda^{-1}}(0,F(v^+)) \;\le\; \norm{P(v^+-v)}_{\Lambda^{-1}}
    \;\le\; \sqrt{\lambda_{\max}}\,\norm{v - T_\omega(v)}_P .
\]
Consequently $\gp(T_\omega(v)) \le \sqrt{\lambda_{\max}}\,\diam_\Lambda(\cZ)\,\norm{v-T_\omega(v)}_P$.
\end{lem}

\begin{proof}
Write $\cS v + \varrho := (c - A^\top y,\ Ax - q)$, so that
$F(v) = \cS v + \varrho + N_{\cZ}(v)$.
Let $u = (x',y') \in \cZ$. Then
\[
    \langle \cS v+\varrho,\ v-u\rangle
    \;=\; \langle c - A^\top y,\ x-x'\rangle + \langle Ax-q,\ y-y'\rangle
    \;=\; L(x,y') - L(x',y),
\]
after cancelling the two copies of $\langle y, Ax\rangle$.
Maximizing over $u \in \cX\times\cY$ decouples into $\max_{\cY}L(x,\cdot)$ and
$-\min_{\cX}L(\cdot,y)$, so
$\max_{u\in\cZ}\langle \cS v+\varrho, v-u\rangle = \gp(v)$.

Let $\xi \in F(v)$ have minimal $\norm{\cdot}_{\Lambda^{-1}}$ norm, say
$\xi = \cS v+\varrho+\zeta$ with $\zeta \in N_{\cZ}(v)$.
For $u \in \cZ$,
\[
    \langle \cS v+\varrho,\ v-u\rangle \;\le\; \langle\xi,v-u\rangle
    \;\le\; \norm{\xi}_{\Lambda^{-1}}\norm{v-u}_{\Lambda}
    \;\le\; \norm{\xi}_{\Lambda^{-1}}\diam_\Lambda(\cZ) .
\]
The first inequality uses $\langle\zeta,u-v\rangle \le 0$.
The second is Cauchy--Schwarz in the $(\Lambda,\Lambda^{-1})$ pairing.
The third uses $u,v\in\cZ$.
Maximizing over $u$ gives the first claim.

For the second, $-P(v^+-v) \in F(v^+)$, because $v^+ = T_\omega(v)$ is the proximal-point step.
From $P \preceq \lambda_{\max}\Lambda$ we get $P\Lambda^{-1}P \preceq \lambda_{\max}P$:
conjugating by $\Lambda^{-1/2}$ turns this into $M^2 \preceq \lambda_{\max}M$ for
$M := \Lambda^{-1/2}P\Lambda^{-1/2}$, which is $M \preceq \lambda_{\max}I$.
Hence $\norm{Pu}_{\Lambda^{-1}}^2 = u^\top P\Lambda^{-1}Pu \le \lambda_{\max}u^\top Pu
= \lambda_{\max}\norm{u}_P^2$.
The last claim combines the two, using $T_\omega(v) \in \cZ$.
\end{proof}

\begin{lem}[Sharpness of the fixed-point residual]\label{lem:sharp}
For every $v \in \RR^{n+m}$,
\[
    \dist_P(v,\cZ^\star) \;\le\; \pran{1 + \gamma\,\diam_\Lambda(\cZ)\,\lambda_{\max}}\,\norm{v-T_\omega(v)}_P .
\]
\end{lem}

\begin{proof}
Write $v^+ := T_\omega(v) \in \cZ$ and let $z^\star$ be the $\Lambda$-nearest point of $\cZ^\star$
to $v^+$. Then
\[
\begin{aligned}
    \dist_P(v,\cZ^\star)
    &\;\le\; \norm{v-v^+}_P + \norm{v^+-z^\star}_P
    \;\le\; \norm{v-v^+}_P + \sqrt{\lambda_{\max}}\,\dist_\Lambda(v^+,\cZ^\star) \\
    &\;\le\; \norm{v-v^+}_P + \sqrt{\lambda_{\max}}\,\gamma\,\gp(v^+)
    \;\le\; \pran{1 + \gamma\,\diam_\Lambda(\cZ)\,\lambda_{\max}}\norm{v-v^+}_P ,
\end{aligned}
\]
where the first inequality is the triangle inequality in $\norm{\cdot}_P$ together with
$z^\star\in\cZ^\star$, the second is $P \preceq \lambda_{\max}\Lambda$ applied to
$v^+-z^\star$, the third is
the definition of $\gamma$ at $v^+ \in \cZ$, and the fourth is the last claim of
Lemma~\ref{lem:gapres}.
\end{proof}

\subsection{Proof of Proposition~\ref{prop:qubo-sharp}}

\begin{proof}[Proof of Proposition~\ref{prop:qubo-sharp}]
The two value polytopes are those of $W$.
The diameter and the canonical norm are those of the instance the method runs on.
Half-integrality is a statement about $W$, and it transports as follows.

The rows cutting out $\cQ_P$ mention no data: covering, linking, $x\in[0,1]^n$ and
$z\in[0,1]^E$.
Normalization therefore fixes $\cQ_P$ and its vertices, and divides its objective by
$\norm{(d;w)}_2$.
The dual value LP is the other way round.
Its objective $\sum_i\varphi_i-\sum_{\Ep}b_e$ has coefficients $\pm1$, so normalization fixes it.
Every datum defining $\cQ_D$ is positively homogeneous of degree one in
$(d,w,a,b,\varphi)$: the bounds $0\le a,b\le w$, the levels $B_i$, and the functions $g_i$.
So $\cQ_D$ becomes $\cQ_D/\norm{(d;w)}_2$, with vertex set that of $W$ divided by
$\norm{(d;w)}_2$.
Either way every vertex value is divided by $\norm{(d;w)}_2$.
Lemmas~\ref{lem:qubo-halfint} and~\ref{lem:qubo-dualmargin} with
Corollary~\ref{cor:lattice} at $q = 2$ then give
$\delta_P,\delta_D \ge 1/(2\norm{(d;w)}_2)$ on the instance the method runs on.

With $\norm{\Lambda_{\cX}}_2 = \max_i\deg(i)$ and $\norm{\Lambda_{\cY}}_2 = 2$,
Lemma~\ref{lem:gammavertex} reads
\begin{equation}\label{eq:gamma-chain}
\begin{aligned}
    \gamma &\;\le\; 2\norm{(d;w)}_2\max\left\{\sqrt{\max_i\deg(i)}\,D_P,\ \sqrt2\,D_D\right\}
      \;\le\; 2\norm{(d;w)}_2\max\left\{\sqrt{2nm},\ \sqrt{10m}\right\} \\
      &\;=\; 2\norm{(d;w)}_2\sqrt{2nm}
      \;\le\; 4\norm{W}_F\sqrt{nm} .
\end{aligned}
\end{equation}
The second inequality bounds the two polytopes.
The equality uses $n \ge 5$.
The last converts $\norm{(d;w)}_2$ to $\norm{W}_F$.

Take the primal polytope first.
$\cQ_P \subseteq [0,1]^{n+m}$ gives $D_P^2 \le n+m$.
The graph is simple on $n$ nodes with $m$ edges, so $\max_i\deg(i) \le \min\{n-1,m\}$.
If $m \ge n$ then $\max_i\deg(i)(n+m) \le (n-1)\cdot 2m \le 2nm$.
If $m < n$ then $\max_i\deg(i)(n+m) \le m\cdot 2n = 2nm$.
Either way $\max_i\deg(i) D_P^2 \le 2nm$.

Now the dual polytope.
$\Em$ and $\Ep$ partition $E$, so $\cQ_D$ sits inside
$\prod_{e\in E}[0,w_e]\times\prod_i[-B_i,0]$.
The method runs on that box divided by $\norm{(d;w)}_2$, so
\[
    D_D^2 \;\le\; \frac{\norm{w}_2^2 + \sum_iB_i^2}{\norm{(d;w)}_2^2}
          \;\le\; \frac{2\norm{d}_2^2 + \pran{1+4\max_i\deg(i)}\norm{w}_2^2}{\norm{(d;w)}_2^2}
          \;\le\; 1 + 4\max_i\deg(i) \;\le\; 5m ,
\]
and $2D_D^2 \le 10m$.
Here $B_i^2 \le 2d_i^2+2\pran{\sum_{e\ni i}w_e}^2$.
Cauchy--Schwarz, counting each edge at both endpoints, gives
$\sum_i\pran{\sum_{e\ni i}w_e}^2 \le \sum_i\deg(i)\sum_{e\ni i}w_e^2
\le 2\max_i\deg(i)\norm{w}_2^2$.
Next, $\max_i\deg(i) \ge 1$ gives $2 \le 1+4\max_i\deg(i)$, which puts
$\norm{(d;w)}_2^2 = \norm{d}_2^2+\norm{w}_2^2$ in the numerator.
Last, $\max_i\deg(i) \le m$ and $m \ge 1$ give $1+4\max_i\deg(i) \le 5m$.
The maximum in \eqref{eq:gamma-chain} is the primal term exactly when $2nm \ge 10m$, that
is when $n \ge 5$.

For the conversion, $d_i = W_{ii}$ and $w_e = 2|W_{ij}|$ give
$\norm{W}_F^2 = \norm{d}_2^2 + \tfrac12\norm{w}_2^2$.
Hence $\norm{(d;w)}_2^2 = \norm{d}_2^2+\norm{w}_2^2 \le 2\norm{W}_F^2$, and
$\norm{(d;w)}_2 \le \sqrt2\norm{W}_F$.

Every node carries an edge, so $\deg(i) \ge 1$ and
\begin{equation}\label{eq:diam-chain}
    \diam_\Lambda(\cZ)^2 \;=\; \sum_i \deg(i) + \frac{2\norm{w}_2^2}{\norm{(d;w)}_2^2}
    \;=\; 2m + \frac{2\norm{w}_2^2}{\norm{(d;w)}_2^2}
    \;\le\; 2m+2 \;\le\; 4m .
\end{equation}
The first equality is the diameter of
$[0,1]^n \times \prod_e\left[0,\ w_e/\norm{(d;w)}_2\right]$.
The second is the handshake identity.
The first inequality is $\norm{w}_2 \le \norm{(d;w)}_2$, and the second is $m \ge 1$.

Combining,
\[
    1 + \gamma\,\diam_\Lambda(\cZ)\,\lambda_{\max}
      \;\le\; 1 + 4\norm{W}_F\sqrt{nm}\cdot 2\sqrt m\cdot\frac2\eta
      \;=\; 1 + \frac{16}{\eta}\norm{W}_F\,m\sqrt n .
\]
The inequality substitutes \eqref{eq:gamma-chain}, \eqref{eq:diam-chain} and
$\lambda_{\max} < 2/\eta$.
With Lemma~\ref{lem:sharp} this is \eqref{eq:sharp-final}.
\end{proof}

\section{Proof of Proposition~\ref{prop:luyang}(b)}\label{app:restart}

Proposition~\ref{prop:luyang}(b) is \cite[Theorem~2]{lu_restarted_2024} carrying the
reflection of \cite[Section~6]{lu_restarted_2024}, stated for $T_\omega$ in the canonical
norm \eqref{eq:canonical}.
Its proof uses three facts: $T_\omega$ is firmly nonexpansive in $\norm{\cdot}_P$ with
fixed point set $\cZ^\star$, Proposition~\ref{prop:luyang}(a), and
Definition~\ref{def:alpha}.
The first holds at every $P \succ 0$ by Section~\ref{sec:alg}, and the second follows
from it by applying \cite{lieder_convergence_2021} to $2T_\omega - I$.
The third is where \cite{lu_restarted_2024} uses the standard form:
\cite[Proposition~2]{lu_restarted_2024} supplies, for each $\rho > 0$, a modulus valid
on $\{v : \norm{T_\omega(v)}_2 \le \rho\}$, which the proof there pairs with a bound on
the iterates.
Definition~\ref{def:alpha} holds at every $v$ and so serves at every $\rho$.
One further difference is in the rule: \eqref{eq:adaptive-restart} applies the residual
test at every epoch, where \cite{lu_restarted_2024} sets the length of the first one in
advance.
The argument below is that of \cite{lu_restarted_2024}, recorded for completeness.

\begin{proof}[Proof of Proposition~\ref{prop:luyang}(b)]
Write $\tau^j$ for the length of epoch $j$, so that \eqref{eq:adaptive-restart} holds at
$k = \tau^j$ and fails at $k = \tau^j - 1$, and $v^{j+1,0} = T_\omega(v^{j,\tau^j})$ by
\eqref{eq:halpern-anchor}.
At any $j$ with $v^{j,0} \in \cZ^\star$ both claims hold with $\tau^j = 0$, so take
$\norm{v^{j,0} - T_\omega(v^{j,0})}_P > 0$, which gives $\tau^j \ge 1$.
Proposition~\ref{prop:luyang}(a) sends $\norm{v^{j,k} - T_\omega(v^{j,k})}_P$ to
zero in $k$, so $\tau^j$ is finite.

For the epoch length,
\[
\begin{aligned}
    \frac1e\norm{v^{j,0} - T_\omega(v^{j,0})}_P
    &\;<\; \norm{v^{j,\tau^j-1} - T_\omega(v^{j,\tau^j-1})}_P
    \;\le\; \frac{\dist_P(v^{j,0},\cZ^\star)}{\tau^j}
    \;\le\; \frac{1}{\alpha\,\tau^j}\norm{v^{j,0} - T_\omega(v^{j,0})}_P ,
\end{aligned}
\]
where the first inequality is the failure of \eqref{eq:adaptive-restart} at
$k = \tau^j - 1$, the second is Proposition~\ref{prop:luyang}(a) at that $k$ with
anchor $v^{j,0}$, and the third is Definition~\ref{def:alpha} at $v^{j,0}$.
Hence $\tau^j < e/\alpha$ and so $\tau^j \le \lceil e/\alpha\rceil$.

Nonexpansiveness of $T_\omega$ and \eqref{eq:adaptive-restart} at $k = \tau^i$ give, for
every $i \ge 0$,
\begin{equation}\label{eq:anchor-residual}
    \norm{v^{i+1,0} - T_\omega(v^{i+1,0})}_P
    = \norm{T_\omega(v^{i,\tau^i}) - T_\omega\pran{T_\omega(v^{i,\tau^i})}}_P
    \le \norm{v^{i,\tau^i} - T_\omega(v^{i,\tau^i})}_P
    \le \frac1e\norm{v^{i,0} - T_\omega(v^{i,0})}_P .
\end{equation}
For the distance,
\[
\begin{aligned}
    \dist_P(v^{j+1,0},\cZ^\star)
    &\;\le\; \frac1\alpha\norm{v^{j+1,0} - T_\omega(v^{j+1,0})}_P
    \;\le\; \frac{e^{-j}}{\alpha}\norm{v^{1,0} - T_\omega(v^{1,0})}_P \\
    &\;\le\; \frac{e^{-j}}{\alpha}\norm{v^{0,\tau^0} - T_\omega(v^{0,\tau^0})}_P
    \;\le\; \frac{e^{-j}}{\alpha\pran{\tau^0+1}}\dist_P(v^{0,0},\cZ^\star) ,
\end{aligned}
\]
where the first inequality is Definition~\ref{def:alpha} at $v^{j+1,0}$, the second
applies \eqref{eq:anchor-residual} at $i = j,\dots,1$, the third is the first two
relations of \eqref{eq:anchor-residual} at $i = 0$, and the last is
Proposition~\ref{prop:luyang}(a) at $k = \tau^0$ with anchor $v^{0,0}$.
\end{proof}

\bibliographystyle{amsplain}
\section{Solver parameters}\label{app:params}

\begin{table}[ht]
    \centering
    \begin{tabular}{l l c c}
        \toprule
        Parameter & Symbol & PRIA (JRP) & PRAQ (QUBO) \\
        \midrule
        Step size & $\eta$ & $0.999/\sqrt{2}$ & $0.999$ \\
        Initial primal weight & $\omega^0$ & $0.108$ & $1.000$ \\
        Sufficient restart & $\beta_{\text{suff}}$ & $e^{-1}$ & $e^{-1}$ \\
        Necessary restart & $\beta_{\text{nec}}$ & $0.571$ & $0.571$ \\
        Artificial restart & $\beta_{\text{art}}$ & $0.247$ & $0.247$ \\
        Initial smoothing exponent & $\theta^0$ & $0.819$ & $0.819$ \\
        Smoothing target & $\theta_{\text{tgt}}$ & $0.428$ & $0.428$ \\
        Smoothing decay & $\theta_{\text{dec}}$ & $0.905$ & $0.905$ \\
        Check cadence & $\ell$ & $16$ & $16$ \\
        \bottomrule
    \end{tabular}
\caption{Default parameters of the method of Section~\ref{sec:alg} for the two instantiations, tuned by randomized search with Optuna \cite{akiba_optuna_2019} on instances disjoint from the benchmark sets.
The step sizes are the two coupling norms of Sections~\ref{sec:qubo} and~\ref{sec:jrp}.
The sufficient restart parameters are the analytic optima for the $\norm{\cdot}_P$ norm of QUBO.
The other parameters are tuned.}
    \label{tab:hyperparameters}
\end{table}

\end{document}